\documentclass[12pt,leqno]{article}

\def\diam{\mathop{\rm diam}}
\def\dist{\mathop{\rm dist}}
\def\Lip{\mathop{\rm Lip}}

\newtheorem{theorem}{Theorem}
\newtheorem{lemma}[theorem]{Lemma}
\newtheorem{proposition}[theorem]{Proposition}
\newtheorem{sublemma}[theorem]{Sublemma}
\newtheorem{definition}[theorem]{Definition}
\newtheorem{corollary}[theorem]{Corollary}
\newtheorem{problem}[theorem]{Problem}
\newtheorem{remark}[theorem]{Remark}
\newtheorem{claim}[theorem]{Claim}
\newtheorem{assumptions}[theorem]{Assumptions}
\newtheorem{examples}[theorem]{Examples}
\newtheorem{basicfact}[theorem]{Basic Fact}

\newcommand{\begintheorem}{\addtocounter{equation}{1}\begin{theorem}}
\newcommand{\beginlemma}{\addtocounter{equation}{1}\begin{lemma}}
\newcommand{\beginproposition}{\addtocounter{equation}{1}\begin{proposition}}
\newcommand{\beginsublemma}{\addtocounter{equation}{1}\begin{sublemma}}
\newcommand{\begindefinition}{\addtocounter{equation}{1}\begin{definition}}
\newcommand{\begincorollary}{\addtocounter{equation}{1}\begin{corollary}}
\newcommand{\beginproblem}{\addtocounter{equation}{1}\begin{problem}}
\newcommand{\beginremark}{\addtocounter{equation}{1}\begin{remark}}
\newcommand{\beginclaim}{\addtocounter{equation}{1}\begin{claim}}
\newcommand{\beginassumptions}{\addtocounter{equation}{1}\begin{assumptions}}
\newcommand{\beginexamples}{\addtocounter{equation}{1}\begin{examples}}
\newcommand{\beginbasicfact}{\addtocounter{equation}{1}\begin{basicfact}}

\begin{document}

\title{Happy fractals and some aspects of analysis on metric spaces}

\author{Stephen Semmes\thanks{This survey was prepared partially in
connection with the trimester ``Heat kernels, random walks, and
analysis on manifolds and graphs'' at the Centre \'Emile Borel,
Institut Henri Poincar\'e, in the Spring of 2002.  This trimester was
organized by P.~Auscher, G.~Besson, T.~Coulhon, and A.~Grigoryan, and
the author was fortunate to be a participant.  The proceedings will be
published in the Contemporary Mathematics series of the American
Mathematical Society, and a report on the trimester can be found in
\cite{SS1}.  Another survey on related themes is \cite{SS2}.  The
author is grateful to an unnamed reader for many helpful comments and
suggestions.}	\\ [10pt]
\textit{dedicated to Leon Ehrenpreis and Mitchel Taibleson}}

\date{}

\maketitle

\begin{abstract}
In this survey we discuss ``happy fractals'', which are complete
metric spaces which are not too big, in the sense of a doubling
condition, and for which there is a path between any two points whose
length is bounded by a constant times the distance between the two
points.  We also review some aspects of basic analysis on metric
spaces, related to Lipschitz functions, approximations and
regularizations of functions, and the notion of ``atoms''.
\end{abstract}

\tableofcontents

\bigskip\bigskip

	As usual, to say that $(M, d(x,y))$ is a metric space means
that $M$ is a nonempty set and that $d(x,y)$ is a nonnegative
real-valued function on $M \times M$ such that $d(x,y) = 0$ if and
only if $x = y$, $d(x,y) = d(y,x)$ for all $x, y \in M$, and
\begin{equation}
	d(x,z) \le d(x,y) + d(y,z)
\end{equation}
for all $x, y, z \in M$ (the triangle inequality).  Here we shall make
the standing assumption that
\begin{equation}
	\hbox{$M$ has at least $2$ elements,}
\end{equation}
to avoid degeneracies.  If $E$ is a nonempty subset
of $M$, then $\diam E$ denotes the \emph{diameter} of $E$, defined by
\begin{equation}
	\diam E = \sup \{d(u,v) : u, v \in E\}.
\end{equation} 
Given $x$ in $M$ and a positive real number $r$, we let $B(x,r)$ and
$\overline{B}(x,r)$ denote the open and closed balls in $M$ with
center $x$ and radius $r$, so that
\begin{equation}
	\enspace
	B(x,r) = \{y \in M : d(x,y) < r\}, \enspace
		\overline{B}(x,r) = \{y \in M : d(x,y) \le r\}.
\end{equation}
Sometimes there might be another metric space $(N, \rho(u,v))$ in
play, and we may introduce a subscript as in $B_N(w,s)$ to indicate in
which metric space the ball is defined.  

	Of course the $n$-dimensional Euclidean space ${\bf R}^n$ with
the standard metric $|x-y|$ is a basic example of a metric space,
which is always good to keep in mind.  Metric spaces associated to
connected graphs will be discussed in Section \ref{graphs}, and the
special case of Cayley graphs from finitely-generated groups will be
reviewed in Section \ref{Finitely-generated groups}.  In Section
\ref{happy fractals} we consider some notions that apply to any metric
space, concerning rectifiable paths in particular.  Sections
\ref{section on Lipschitz retracts} -- \ref{Some happy fractals from
Helsinki} deal with related notions and examples.  The remaining sections
deal with various general aspects of analysis on metric spaces.

\section{Graphs}
\label{graphs}
\setcounter{equation}{0}

	Suppose that we have a graph consisting of a nonempty set $V$
of vertices and a set $E$ of edges.  An element of $E$ can be
described by an unordered pair of distinct elements of $V$; we do not
wish to consider edges which form loops by themselves, or multiple
edges between the same pair of vertices.  Two vertices connected by an
edge are said to be \emph{adjacent}.

	Let us assume that our graph is \emph{connected}, which is to
say that every pair of vertices can be connected by a finite path.
The length of a path is defined to be the number of edges that the
path traverses.  Thus the length of a path is a nonnegative integer,
which is $0$ in the case of a path that consists of a single vertex
and traverses no edges.

	We define a metric $d(v, w)$ on $V$ by taking $d(v, w)$ to be the
length of the shortest path between $v$ and $w$.  It is easy to see that 
$(V, d(v, w))$ is indeed then a metric space.

	Let us also assume that the graph is \emph{locally finite}, which
is to say that there are only finitely many vertices adjacent to a given
vertex.  For each $p$ in $V$ and each positive integer $m$ one can show
that there are only finitely many vertices whose distance to $p$ is at
most $m$.

	In fact, let us assume that there is a nonnegative integer $k$
such that for every vertex $v$ in $V$ there are at most $k$ vertices
$w$ in $V$ which are adjacent to $v$.  If $k = 0$ then $V$ contains
only one vertex and there are no edges, and if $k = 1$ then $V$ has
either one or two vertices, with no edges if there is only one vertex
and exactly one edge when there are two vertices.  For simplicity let
us assume that $k \ge 2$.

	If $p$ is an element of $V$ and $m$ is a nonnegative integer,
then we define $A_m(p)$ to be the number of vertices $v$ in $V$ whose
distance to $p$ is exactly equal to $m$.  Thus $A_0(p) = 1$, since
$p$ is the only vertex at distance $0$ from itself, and $A_1(p) \le k$,
since $A_1(p)$ is the same as the number of vertices in $V$ which are
adjacent to $p$.  For $m \ge 2$ we have that
\begin{equation}
	A_m(p) \le (k-1) \, A_{m-1}(p).
\end{equation}
Indeed, suppose that $v$ is an element of $V$ whose distance to $p$ is
exactly equal to $m$.  Then there is vertex $w$ in $V$ such that $v$
is adjacent to $w$ and the distance from $w$ to $p$ is exactly $m -
1$.  Since $m \ge 2$, there is also a vertex $u$ in $V$ such that $w$
is adjacent to $u$ and the distance from $u$ to $p$ is exactly $m -
2$.  The total number of vertices in $V$ which are adjacent to $w$ and
which have distance to $p$ equal to $m$ is at most $k - 1$, because
there are at most $k$ vertices which are adjacent to $w$ at all, and
$u$ is adjacent to $w$ and has distance to $p$ equal to $m - 2$.
There are $A_{m-1}(p)$ vertices $w$ whose distance to $p$ is equal to
$m-1$, and hence there are at most $(k - 1) \, A_{m-1}(p)$ vertices
whose distance to $p$ is equal to $m$, as desired.  

	Thus $A_m(p)$ grows at most exponentially in $m$ in general,
and exponential growth is certainly possible, at least when $k \ge 3$.
Of course there are many interesting situations where the growth is in
fact bounded by a polynomial.  In this survey we shall focus on
situations with polynomial growth, and the \emph{doubling condition}
described in Section \ref{happy fractals} gives a nice version of this
which makes sense in any metric space.

	Instead of looking at rates of growth in terms of $A_m(p)$,
one also frequently considers the quantity
\begin{equation}
	\sum_{j=0}^m A_j(p),
\end{equation}
which is the same as the number of elements of $V$ whose distance to
$p$ is at most equal to $m$.

	As a basic example, fix a positive integer $n$, and consider
the set ${\bf Z}^n$ of points in ${\bf R}^n$ with integer coordinates
as a set of vertices.  Two points $v$, $w$ in ${\bf Z}^n$ can be
defined to be adjaced if $v - w$ has $n - 1$ coordinates equal to $0$
and the remaining coordinate equal to $\pm 1$.  This is the same as
saying that $v$, $w$ are adjacent if and only if $|v - w| = 1$.  In this
case it is not difficult to determine the metric on ${\bf Z}^n$ coming
from paths in the graph, namely
\begin{equation}
	d(v, w) = \sum_{j=1}^n |v_j - w_j|,
\end{equation}
where $v_j$, $w_j$ denote the $j$th coordinates of $v$, $w$,
respectively.  This is often called the \emph{taxicab metric},
and it satisfies the following comparison with the Euclidean
distance:
\begin{equation}
	|v - w| \le d(v, w) \le \sqrt{n} \, |v-w|.
\end{equation}
The first inequality can be derived from the triangle inequality for
the standard distance, since each step of size $1$ in the graph metric
is also a step of size $1$ in the Euclidean distance.  The second
inequality is a consequence of the Cauchy--Schwarz inequality.

	In this case the growth is polynomial, with the number of
points at distance to a fixed point $p$ less than or equal to $r$ is
on the order of $r^n$.  Notice that this number does not depend on
$p$, because of translation-invariance.

	Concerning analysis and geometry of graphs and related
matters, see \cite{Thierry, PW, Woess}.

\section{Finitely-generated groups}
\label{Finitely-generated groups}
\setcounter{equation}{0}

	A very interesting special case of graphs and their geometry
comes from \emph{Cayley graphs} of \emph{finitely generated groups}.
Let $\Gamma$ be a group with a finite set $F$ of generators.  Thus
every element of $\Gamma$ can be expressed as a product of elements of
$F$ and their inverses, with the identity element viewed as an
empty product of generators.  For the Cayley graph of $\Gamma$
we use $\Gamma$ as the set of vertices, and define two elements
$\gamma_1$, $\gamma_2$ of $\Gamma$ to be adjacent if one of them
can be written as the product of the other times an element of $F$,
where the group operation is applied in that order.  From this it
follows that the graph is invariant under left-translations, which
is to say that $\gamma_1$, $\gamma_2$ are adjacent if and only if
$\alpha \, \gamma_1$, $\alpha \, \gamma_2$ for any $\alpha$ in $\Gamma$.

	Every pair of elements of $\Gamma$ can be joined by a path
in the Cayley graph, because of the assumption that every element
of $\Gamma$ can be expressed as a product of generators and their
inverses.  If $d(\gamma_1, \gamma_2)$ denotes the distance function
on $\Gamma$ coming from the Cayley graph, then we have that
\begin{equation}
	d(\alpha \, \gamma_1, \alpha \, \gamma_2) = d(\gamma_1, \gamma_2)
\end{equation}
for all $\alpha$, $\gamma_1$, and $\gamma_2$ in $\Gamma$, by
left-invariance of the Cayley graph.

	If $\gamma$ is an element of $\Gamma$, then the number of
elements of $\Gamma$ which are adjacent to $\gamma$ is at most twice
the number of elements of $F$, by construction.  As in the preceding
section, this leads to a simple exponential bound on the growth of the
Cayley graph of $\Gamma$.  Exponential growth occurs for free groups
with at least two generators, and more generally for nonelementary
hyperbolic groups in the sense of Gromov, as in \cite{CDP, CP, GH,
Gromov1, Gromov3}.  Hyperbolic groups have very interesting spaces at
infinity associated to them which satisfy the doubling property
described in the next section.  In addition to the references already
mentioned, see \cite{Coornaert, GP, Pierre1, Pierre2} in this regard.
Note that fundamental groups of compact Riemannian manifolds without
boundary and with strictly negative sectional curvatures are
nonelementary hyperbolic groups.  Simply-connected symmetric spaces
always have compact quotients by a well-known result of Borel
\cite{Borel1, Raghunathan}, and for symmetric spaces of noncompact
type and rank $1$ the sectional curvatures are strictly negative.

	The graph associated to ${\bf Z}^n$ in the previous section
is exactly its Cayley graph as group with the $n$ standard generators,
where each generator has one coordinate equal to $1$ and the others
equal to $0$.  This graph has polynomial growth, as we saw, and
more generally it is a well-known result that the Cayley graph of a
finitely-generated group has polynomial growth when the group is
\emph{virtually nilpotent}, which means that the group contains
a nilpotent subgroup of finite index.  A famous theorem of Gromov
\cite{Gromov2} states that the converse is true.

\section{Happy fractals}
\label{happy fractals}
\setcounter{equation}{0}

	Let us say that a metric space $(M, d(x,y))$ is a \emph{happy
fractal} if the following three conditions are satisfied.  First, $M$
is complete as a metric space.  Second, there is a constant $C_1 > 0$
so that for each pair of points $x$, $y$ in $M$ there is a path in $M$
connecting $x$ to $y$ with length at most $C_1 \, d(x,y)$.  Third, $M$
satisfies the \emph{doubling property} that there is a constant $C_2$
so that any ball $B$ in $M$ can be covered by a family of balls with
half the radius of $B$ and at most $C_2$ elements.

	One might prefer the name \emph{happy metric space}, since the
metric space need not be fractal, as in the case of ordinary Euclidean
spaces.  There are plenty of examples which are more intricate and not
fractal, such as domains or surfaces with cusps.  There can be
interesting fractal behavior at some kind of boundary, if not for the
space itself.

	Cantor sets and snowflake curves give examples of self-similar
fractals which satisfy the doubling condition but are not happy fractals,
because every curve of finite length in these spaces is constant.
Some basic examples of happy fractals will be discussed in the next
few sections.

	It does not seem to be known whether every compact connected
$4$-dimensional topological manifold can be realized as a happy
fractal, i.e., whether every compact Hausdorff topological space which
is locally homeomorphic to the open unit ball in ${\bf R}^4$ has a
topologically-equivalent metric in which it becomes a happy fractal.
This is true for dimensions not equal to $4$, since $n$-dimensional
topological manifolds admit unique smooth structures when $n \le 3$
and they admit unique Lipschitz structures when $n \ge 5$.  See
\cite{Bing, Do-Su, F-Q, Moise, Sullivan} concerning these topics.

	In general dimensions there are plenty of questions about
noncompact spaces.  For instance, in this connection one might
consider conditions of bounded local geometry, with the happy fractal
aspect being concerned with larger scales.  In dimension $4$, let us
recall a well-known result of Quinn that every connected
$4$-dimensional topological manifold can be smoothed in the complement
of a single point.  Of course, near that point there can be a lot of
complications, although there are also topological restrictions since
that point is a topological manifold point.

	To be more precise, a \emph{path} in $M$ which goes from a
point $x$ to a point $y$ is a continuous mapping $p(t)$ defined on
a closed interval $[a,b]$ in the real line and with values in $M$
such that $p(a) = x$ and $p(b) = y$.  If
\begin{equation}
	a = t_0 < t_1 < t_2 < \cdots < t_m = b
\end{equation}
is a partition of $[a,b]$, then we can associate to this partition the
quantity
\begin{equation}
\label{approx to length along partition}
	\sum_{j=1}^m d(p(t_j), p(t_{j-1})),
\end{equation}
which is the approximation to the length of $p$ corresponding to this
partition.  The \emph{length} of the path is defined to be the
supremum of (\ref{approx to length along partition}) over all
partitions of $[a,b]$.  In general this can be infinite.

	A standard observation is that the quantity (\ref{approx to
length along partition}) can only increase as points are added to
the partition, because of the triangle inequality.  Any two partitions
admit a common refinement, for which the approximation to the length
is then greater than or equal to the approximations to the length
associated to the original refinements.

	Suppose that the length of the path $p(t)$ is finite.  Then
the length of the restriction of $p$ to any subinterval of $[a,b]$ is
also finite, and is less than or equal to the length of the whole
path.  Let us define a function $L(u,v)$ for $u, v \in [a,b]$, $u \le
v$, to be the length of the restriction of $p(t)$ to $[u,v]$.  Of
course a constant path has length $0$, which includes the case where
the domain has one element.  Note that
\begin{equation}
\label{d(p(u),p(v)) le L(u,v)}
	d(p(u),p(v)) \le L(u,v)
\end{equation}
for all $u, v \in [a,b]$ with $u \le v$.  If
\begin{equation}
	a \le u \le v \le w \le b,
\end{equation}
then it is not hard to verify that
\begin{equation}
	L(u,w) = L(u,v) + L(v,w),
\end{equation}
using the monotonicity properties of the length, and the possibility
of taking refinements of the partitions in particular.

	Fix $t \in [a,b]$.  If $t > a$, then 
\begin{equation}
	\lim_{s \to t-} L(s,t) = 0.
\end{equation}
This is equivalent to saying that
\begin{equation}
	\lim_{s \to t-} L(a,s) = L(a,t).
\end{equation}
From the definition we know that $L(a,s)$ is monotone increasing in
$s$, so that the limit on the left side exists and is less than or
equal to the right side.  To show that equality holds, one can choose
a partition of $[a,t]$ so that the approximation to the length of
$p(u)$ along this partition is close to $L(a,t)$, and then check that
$L(a,s)$ is greater than or equal to this approximation minus a small
number when $s$ is sufficiently close to $t$.  This employs the
continuity of $p(u)$ at $t$, to move the last point in the partition
from $t$ to $s$ without making more than a small change to the
approximation to the length.

	If $t < b$, then 
\begin{equation}
	\lim_{s \to t+} L(t,s) = 0.
\end{equation}
This is equivalent to
\begin{equation}
	\lim_{s \to t+} L(s,b) = L(t,b),
\end{equation}
which can be verified in the same manner as before.

	Set $\lambda = L(a,b)$, and consider the real-valued function
$\sigma(t)$ defined on $[a,b]$ by
\begin{equation}
	\sigma(t) = L(a,t).
\end{equation}
Thus $\sigma(t)$ is monotone increasing (and not necessarily strictly
increasing), $\sigma(0) = 0$, $\sigma(b) = \lambda$, and $\sigma(t)$
is continuous by the preceding remarks.

	There is a mapping $\widetilde{p} : [0,\lambda] \to M$
such that
\begin{equation}
\label{widetilde{p}(sigma(t)) = p(t)}
	\widetilde{p}(\sigma(t)) = p(t)
\end{equation}
for all $t \in [a,b]$.  In other words, if $s, t \in [a,b]$, $s < t$,
and $\sigma(s) = \sigma(t)$, then $L(s,t) = 0$, so that $p$ is
constant along $[s,t]$, and (\ref{widetilde{p}(sigma(t)) = p(t)})
leads to a single value for $p$ at $\sigma(s) = \sigma(t)$.  Moreover,
(\ref{d(p(u),p(v)) le L(u,v)}) implies that
\begin{equation}
	d(\widetilde{p}(r),\widetilde{p}(w)) \le |r-w|
\end{equation}
for all $r, w \in [0,\lambda]$.  

	On the other hand, if $q : [c,d] \to M$ is a path
such that
\begin{equation}
\label{d(q(s),q(t)) le k |s-t|}
	d(q(s),q(t)) \le k \, |s-t|
\end{equation}
for some constant $k$ and all $s, t \in [c,d]$, then it is easy
to check that the length of $q$ on $[c,d]$ is at most $k \, |c-d|$.
One can trade between $k$ and $|c-d|$ by rescaling in the domain.

	Thus there is a path in $M$ from $x$ to $y$ with length less
than or equal to a constant $A$ if and only if there is a mapping $q :
[0,1] \to M$ such that $q(0) = x$, $q(1) = y$, and (\ref{d(q(s),q(t))
le k |s-t|}) holds for all $s, t \in [0,1]$ with $k \le A$.  Assuming
that there is a path in $M$ from $x$ to $y$ with finite length and
that closed and bounded subsets of $M$ are compact, one can use the
Arzela--Ascoli theorem to find such a mapping $q$ with $k$ as small as
possible, and this minimal $k$ is the same as the length of the
shortest path in $M$ from $x$ to $y$.

	A well-known result in basic analysis states that if $(M,
d(x,y))$ is a complete metric space, then a closed subset $K$ of $M$
is compact if and only if $K$ is totally bounded, which means that for
every $\epsilon > 0$ there is a finite family of balls in $M$ with
radius $\epsilon$ whose union contains $K$.  Thus, if $(M, d(x,y))$ is
complete, then closed and bounded subsets of $M$ are compact if and
only if all balls in $M$ are totally bounded.  It is easy to verify
that the latter holds when $M$ satisfies the doubling property.  In
short, closed and bounded sets are compact in a happy fractal (or
happy metric space).

\section{Lipschitz retracts}
\label{section on Lipschitz retracts}
\setcounter{equation}{0}

	Suppose that $(M, d(x,y))$ is a metric space, and that $A$ and
$E$ are subsets of $M$, with $E \subseteq A$.  A mapping $\phi : A \to
E$ is said to be a \emph{Lipschitz retract} of $A$ onto $E$ if
\begin{equation}
	\phi(x) = x \quad\hbox{for all } x \in E
\end{equation}
and $\phi$ is Lipschitz, so that there is a constant $k \ge 0$ such
that
\begin{equation}
	d(\phi(y), \phi(z)) \le k \, d(y,z)
\end{equation}
for all $y, z \in A$.  Note that if $M$ is complete and $E$ is a
closed subset of $M$, then one can always take $A$ to be closed,
because any Lipschitz mapping from $A$ into $E$ can be extended to a
Lipschitz mapping from the closure of $A$ into $E$, and with the same
Lipschitz constant $k$.

	Let us say that a complete metric space $(N, \rho(u,v))$ is a
\emph{Lipschitz extension space} with constant $s \ge 1$ if for every
separable metric space $(M, d(x,y))$ and every mapping $f$ from a
subset $Z$ of $M$ into $N$ which is Lipschitz with constant $L$, so
that
\begin{equation}
	\rho(f(x), f(y)) \le L \, d(x,y)
\end{equation}
for all $x, y \in Z$, there is an extension of $f$ to a Lipschitz
mapping from $M$ into $N$ with constant $s \, L$.  

\beginremark
{\rm If $(M, d(x,y))$ is a separable metric space and $E$ is a
subset of $M$, and if $(E, d(x,y))$ satisfies the Lipschitz extension
property with constant $s$, then there is a Lipschitz retraction from
$M$ onto $E$ with constant $s$, simply by extending the identity mapping
on $E$. }
\end{remark}

	The requirement above that $N$ be complete is not really
needed, since it can be derived from the extension property.  The
restriction to metric spaces $M$ which are \emph{separable} --- i.e.,
which contain a countable dense subset --- is made because we shall
only be concerned with spaces that satisfy this condition, and because
it permits one to avoid such things as transfinite induction.
Specifically, one can make the following observation.

\beginlemma
\label{reformulation of Lipschitz extension property}
Let $(N, \rho(u,v))$ be a complete metric space.  A necessary and
sufficient condition for $N$ to satisfy the Lipschitz extension
property with constant $s$ is that it satisfy this property in the
special case where the metric space $(M, d(x,y))$ and the subset $Z$
of $M$ have the feature that $M \backslash Z$ is at most countable.
\end{lemma}

	Indeed, given arbitrary $(M, d(x,y))$, $Z$, $f$, and $L$ as in
the definition of the Lipschitz extension property, one can first use
separability of $M$ to find a subset $M_0$ of $M$ such that $M_0$
contains $Z$, $M_0 \backslash Z$ is at most countable, and $M_0$ is
dense in $M$.  Under the restricted version of the Lipschitz extension
property mentioned in the lemma, one can extend $f$ to a Lipschitz
mapping from $M_0$ to $N$ with Lipschitz constant $s \, L$.  The
completeness of $N$ then permits this mapping to be extended to one
from all of $M$ into $N$, with Lipschitz constant $s \, L$ still.

	I learned the next lemma from M.~Gromov, as well as the way it
can be used.

\beginlemma
\label{criterion for Lipschitz extension property with s = 1}
Let $(N, \rho(u,v))$ be a complete metric space.  A necessary and
sufficient condition for $(N, d(x,y))$ to satisfy the Lipschitz
extension property with constant $s = 1$ is that it satisfy this
property in the special case where the metric space $(M, d(x,y))$
and the subset $Z$ of $M$ have the feature that $M \backslash Z$
contains only one element.
\end{lemma}

	Indeed, if one can extend a Lipschitz mapping to a set with
one extra element, without increasing the Lipschitz constant, then one
can repeat this to get extensions to sets with arbitrary finite
numbers of additional elements, or even countably many additional
elements, without increasing the Lipschitz constant.  The preceding
lemma then applies to deal with the general case.

\beginlemma
\label{criterion for the previous criterion}
Let $(N, \rho(u,v))$ be a complete metric space.
Suppose that for every collection
\begin{equation}
	\{B_i\}_{i \in I} = \{\overline{B}_N(u_i, r_i)\}_{i \in I}
\end{equation}
of closed balls in $N$ such that $I$ is at most countable and
\begin{equation}
\label{rho(u_i, u_j) le r_i + r_j hbox{for all } i, j in I}
	\rho(u_i, u_j) \le r_i + r_j \quad\hbox{for all } i, j \in I
\end{equation}
we have that
\begin{equation}
	\bigcap_{k \in I} B_k \ne \emptyset.
\end{equation}
Then $(N, \rho(u,v))$ satisfies the Lipschitz extension property with
$s = 1$.
\end{lemma}

	Note that the completeness of $N$ corresponds in fact to the
special case of the condition in the lemma where $\{B_i\}_{i \in I}$
is a sequence of closed balls which is decreasing in terms of
inclusion and whose radii are tending to $0$.

	To prove the lemma, it is enough to obtain one-point
extensions, as in Lemma \ref{criterion for Lipschitz extension
property with s = 1}.  Let $(M, d(x,y))$, $Z$, $f$, and $L$
be given as in the definition of the Lipschitz extension property,
with $M \backslash Z$ containing exactly one element $w$.  For 
each $z \in Z$, consider the closed ball 
\begin{equation}
	B_z = \overline{B}_N(f(z), L \, d(w,z))
\end{equation}
in $N$.  If $z_1, z_2 \in Z$, then
\begin{equation}
	\rho(f(z_1), f(z_2)) \le L \, d(z_1, z_2) 
			\le L \, d(w,z_1) + L \, d(w,z_2).
\end{equation}
In other words, this family of balls satisfies the condition
(\ref{rho(u_i, u_j) le r_i + r_j hbox{for all } i, j in I}) in Lemma
\ref{criterion for the previous criterion}.  Although $Z$ may not be at
most countable, one can use the separability of $M$ to obtain that there
is a dense subset $I$ of $Z$ which is at most countable.  The hypothesis
of the lemma then implies that
\begin{equation}
	\bigcap_{z \in I} B_z \ne \emptyset.
\end{equation}
Fix a point $\alpha$ in this intersection, and set $f(w) = \alpha$.
We have that
\begin{equation}
\label{rho(f(w), f(z)) = rho(alpha, f(z)) le L d(w,z)}
	\rho(f(w), f(z)) = \rho(\alpha, f(z)) \le L \, d(w,z)
\end{equation}
for all $z in I$, precisely because $\alpha \in B_z$ for all $z \in
I$.  By continuity, (\ref{rho(f(w), f(z)) = rho(alpha, f(z)) le L
d(w,z)}) holds for all $z \in Z$.  Thus we have an extension of $f$ to
$M = Z \cup \{w\}$ which is Lipschitz with constant $L$, as desired.
This proves the lemma.

\begincorollary
\label{Lipschitz extension property for {bf R} with s = 1}
The real line ${\bf R}$ with the standard metric $|x-y|$ satisfies
the Lipschitz extension property with $s = 1$.
\end{corollary}

	Of course this is well-known and can be established by other
means, as in Section \ref{More on Lipschitz functions}, but one can
check that the hypothesis of Lemma \ref{criterion for the previous
criterion} holds in this case.  To be more precise, the $B_i$'s are
closed and bounded intervals in this case, and the condition
(\ref{rho(u_i, u_j) le r_i + r_j hbox{for all } i, j in I}) implies
that every pair of these intervals intersects.  The special geometry
of the real line implies that the intersection of all of the intervals
is nonempty.

	Part of the point of this kind of approach is that it can be
applied to tree-like spaces.  As a basic scenario, suppose that $(T,
\sigma(p,q))$ is a metric space which consists of a finite number of
pieces which we shall call \emph{segments}, and which are individually
isometrically equivalent to a closed and bounded interval in the real
line.  We assume that any two of these segments are either disjoint or
that their intersection consists of a single point which is an endpoint
of each of the two segments.  We also ask that $T$ be connected, and
that the distance between any two elements of $T$ is the length of
the shortest path that connects them.  One may as well restrict one's
attention to paths which are piecewise linear, and the length of the
paths is easy to determine using the fact that each segment is
equivalent to a standard interval (of some length).

	So far these conditions amount to saying that $T$ is a finite
graph, with the internal geodesic distance.  Now let us also ask that
$T$ be a \emph{tree}, in the sense that any simple closed path in $T$
is trivial, i.e., consists only of a single point.

	The effect of this is that if $p$ and $q$ are elements of $T$,
then there is a special subset $S(p,q)$ of $T$ which is
isometrically-equivalent to a closed and bounded interval in the real
line, with $p$ and $q$ corresponding to the endpoints of this
interval.  In practice, with a typical picture of a tree, it is very
easy to draw the set $S(p,q)$ for any choice of $p$ and $q$.  This set
gives the path of minimal length between $p$ and $q$ (through the
isometric equivalence mentioned before), and it satisfies a stronger
minimality property, namely, any path in $T$ connecting $p$ and $q$
contains $S(p,q)$ in its image.

\beginlemma
\label{the case of finite trees}
Under the conditions just described, $(T, \sigma(p,q))$ satisfies
the hypotheses of Lemma \ref{criterion for the previous criterion}.
As a result, $(T, \sigma(p,q))$ enjoys the Lipschitz extension property
with $s = 1$.
\end{lemma}

	Clearly $T$ is complete, and in fact compact.  Now suppose
that $\{B_i\}_{i \in I}$ is a family of closed balls in $T$.  The
condition (\ref{rho(u_i, u_j) le r_i + r_j hbox{for all } i, j in I})
implies in this setting (and in any geodesic metric space) that
any two of the $B_i$'s intersect.  (Note that the converse always
holds.)

	Let us call a subset $C$ of $T$ \emph{convex} if $p, q \in C$
implies that $S(p,q) \subseteq C$.  Of course convexity in this sense
implies connectedness, and in fact connectedness implies convexity
because of the assumption that $T$ is a tree.  That is, $S(p,q)$ is
contained in any connected set that contains $p$ and $q$.  Of course
connected subsets of $T$ have a simple structure, since a connected
subset of an interval in the real line is also an interval.

	Because the distance on $T$ is defined in terms of lengths of
paths, open and closed balls in $T$ are connected, and hence convex.
The intersection of two convex sets is also convex, by definition.

	Suppose that $C_1$, $C_2$, and $C_3$ are convex subsets of $T$
such that $C_1 \cap C_2$, $C_1 \cap C_3$, and $C_2 \cap C_3$ are all
nonempty.  Let us check that $C_1 \cap C_2 \cap C_3$ is nonempty as
well.  Let $p_{12}$, $p_{13}$, and $p_{23}$ be elements of $C_1 \cap
C_2$, $C_1 \cap C_3$, and $C_2 \cap C_3$, respectively.  Observe that
\begin{equation}
\label{S(p_{12}, p_{13}) subseteq S(p_{12}, p_{23}) cup S(p_{23}, p_{13})}
	S(p_{12}, p_{13}) \subseteq S(p_{12}, p_{23}) \cup S(p_{23}, p_{13}),
\end{equation}
since the right side defines a connected subset of $T$ that contains
$p_{12}$ and $p_{13}$.  As before, there is an isometric equivalence
between $S(p_{12}, p_{13})$ and a closed and bounded interval $I$ in
the real line, where $p_{12}$, $p_{13}$ correspond to the endpoints of
$I$.  On the other hand, $S(p_{12}, p_{13}) \cap S(p_{12}, p_{23})$
and $S(p_{12}, p_{13}) \cap S(p_{23}, p_{13})$ are closed convex
subsets of $S(p_{12}, p_{13})$, and hence correspond to closed
subintervals $J$, $K$ of $I$.  From (\ref{S(p_{12}, p_{13}) subseteq
S(p_{12}, p_{23}) cup S(p_{23}, p_{13})}) we obtain that $I \subseteq
J \cup K$, which implies that $J \cap K \ne \emptyset$, since $J$ and
$K$ are closed.  Any element of $J \cap K$ corresponds to a point
in $S(p_{12}, p_{13})$ that also lies in $S(p_{12}, p_{23})$ and
$S(p_{23}, p_{13})$.  Because $C_1$, $C_2$, and $C_3$ are convex,
$S(p_{12}, p_{13}) \subseteq C_1$, $S(p_{12}, p_{23}) \subseteq C_2$,
and $S(p_{23}, p_{13}) \subseteq C_3$.  In other words, we get an
element of the intersection of $C_1$, $C_2$, and $C_3$, as desired.

	Because the intersection of convex sets is convex, one can
iterate this result to obtain that if $C_1, C_2, \ldots, C_\ell$ are
convex sets in $T$ such that the intersection of any two of them is
nonempty, then $\bigcap_{i=1}^\ell C_i \ne \emptyset$.  For closed
convex sets, which are then compact since $T$ is compact, one can get
the same result for an infinite family of convex sets.  This uses the
well-known general result that the intersection of a family of compact
sets is nonempty if the intersection of every finite subfamily is
nonempty.

	This shows that $(T, \sigma(p,q))$ satisfies the hypotheses of
Lemma \ref{criterion for the previous criterion}, since closed balls
are closed convex sets.  This completes the proof of Lemma \ref{the
case of finite trees}.

	Of course there are analogous results for more complicated
trees or tree-like sets.  Let us note that one might have the set
sitting inside of a Euclidean space, but where the internal geodesic
metric is not quite the same as the restriction of the ambient
Euclidean metric.  If the two are comparable, in the sense that each
is bounded by a constant multiple of the other, then the Lipschitz
extension property for one metric follows from the same property for
the other metric, with a modestly different constant.

\section{The Sierpinski gasket and carpet}
\label{The Sierpinski gasket and carpet}
\setcounter{equation}{0}

	The \emph{Sierpinski gasket} is the compact set in ${\bf R}^2$
which is constructed as follows.  One starts with the unit equilateral
triangle, with bottom left vertex at the origin and bottom side along
the $x_1$-axis.  By ``triangle'' we mean the closed set which includes
both the familiar polygonal curve and its interior.  This triangle can
be subdivided into four parts each with sidelength equal to half of
the original.  The vertices of the four new triangles are vertices of
the original triangle or midpoints of its sides.  One removes the interior
of the middle triangle, and keeps the other three triangles in the first
stage.  One then repeats the process for each of those triangles, and
so on.  The Sierpinski gasket is the compact set without interior
which occurs in the limit, and which is the intersection of the sets
which are finite unions of triangles which occur at the finite stages
of the construction.

	Similarly, the \emph{Sierpinski carpet} is the compact set in
${\bf R}^2$ defined in the following manner.  One starts with the unit
square, where ``square'' also means the familiar polygonal curve
together with its interior.  One decomposes the unit square into
nine smaller squares, each with sidelength equal to one-third that
of the original.  One removes the interior of the middle square,
and keeps the remaining eight squares for the first stage of the
construction.  One then repeats the construction for each of the
smaller squares, and so on.  The Sierpinski carpet is the compact
set without interior which occurs in the limit and is the intersection
of the sets which are the finite unions of squares from the finite
stages of the construction.

	The Sierpinski gasket and carpet provide well-known basic
examples of happy fractals.  The main point is that if $x$, $y$ are
two elements of one of these sets, then $x$ and $y$ can be connected
by a curve in the set whose length is bounded by a constant times
$|x-y|$.  This is not too difficult to show, using the sides of the
triangles and squares to move around in the sets.  

	For neither of these sets is there a continuous retraction
(let alone a Lipschitz retraction) from ${\bf R}^2$ onto the set.
There is not even a continuous retraction from a neighborhood of the
set onto the set.  This is because in both cases there are arbitrarily
small topological loops, given by boundaries of triangles or squares,
which cannot be contracted to a point in the set, but can easily be
contracted to a point in ${\bf R}^2$, within the particular triangle
or square.  If there were a retraction whose domain included suc a
triangle or square, then the contraction of the loop could be pushed
back into the Sierpinski gasket or carpet, where in fact it cannot
exist.

	However, one can retract the complement of a triangle or
square onto its boundary.  If one removes a hole from each open
triangle or square in the complement of the Sierpinski gasket or
carpet, then one can define a continuous retraction on the fatter sets
that one obtains, i.e., as the complement of the union of the holes.
The domain of the retraction is reasonably fat, but it still does not
contain a neighborhood of the Sierpinski gasket or carpet.  If one is
careful to choose the holes so that they always contain a disk of
radius which is greater than or equal to a fixed positive constant
times the diameter of the corresponding triangle or square, then one
can get a Lipschitz retraction.

	There are also nice Lipschitz retractions from the Sierpinski
gasket or carpet onto subsets of itself.  For instance, one can start
by pushing parts of the gasket or carpet in individual triangles or
squares to all or parts of the boundaries of these triangles or
squares.  One can often move what remains into the rest of the gasket
or carpet that is not being moved.

\section{Heisenberg groups}
\label{Heisenberg groups}
\setcounter{equation}{0}

	Let $n$ be a positive integer.  Define $H_n$ first as a set by
taking ${\bf C}^n \times {\bf R}$, where ${\bf C}$ denotes the complex
numbers.  The group law is given by
\begin{equation}
	(w,s) \circ (z,t) = 
   \biggl(w + z, s + t + 2 \, {\rm Im} \sum_{j=1}^n w_j \overline{z_j} \biggr),
\end{equation}
where ${\rm Im} \, a$ denotes the imaginary part of a complex number
$a$, and $w_j$, $z_j$ denote the $j$th components of $w, z \in {\bf
C}^n$.

	It is not difficult to verify that this does indeed define a
group structure on $H_n$.  In this regard, notice that the inverse of
$(w,s)$ in $H_n$ is given by 
\begin{equation}
	(w,s)^{-1} = (-w, -s).
\end{equation}

	For each positive real number $r$, define the ``dilation''
$\delta_r$ on $H_n$ by
\begin{equation}
	\delta_r(w,s) = (r \, w, r^2 \, s).
\end{equation}
One can check that these dilations define group automorphisms of $H_n$,
i.e.,
\begin{equation}
	\delta_r\bigl((w,s) \circ (z,t) \bigr) =
		\delta_r(w,s) \circ \delta_r(z,t).
\end{equation}
Also, for $r_1, r_2 > 0$ we have that
\begin{equation}
	\delta_{r_1}(\delta_{r_2}(w,s)) = \delta_{r_1 \, r_2}(w,s).
\end{equation}
Let us note that the group law and the dilations are compatible with
the standard Euclidean topology on $H_n$, i.e., they define continuous
mappings.

	Let us call a nonnegative real-valued function $N(\cdot)$ on
$H_n$ a \emph{norm} if it satisfies the following conditions: (a) $N$
is continuous; (b) $N$ takes the value $0$ at the origin and is
strictly positive at other points in $H_n$; (c)
$N\bigl((w,s)^{-1}\bigr) = N(w,s)$ for all $(w,s) \in H_n$; (d)
$N(\delta_r(w,s)) = r \, N(w,s)$ for all $r > 0$ and $(w,s) \in H_n$;
and (e) $N$ satisfies the triangle inequality with respect to the
group structure on $H_n$, which is to say that
\begin{equation}
\label{triangle inequality for norm on H_n}
	N\bigl((w,s) \circ (z,t)\bigr) \le N(w,s) + N(z,t)
\end{equation}
for all $(w,s), (z,t) \in H_n$.

	In many situations it is sufficient to work with a weaker
notion, in which (\ref{triangle inequality for norm on H_n}) is
replaced by the ``quasitriangle inequality'' which says that there is
a positive constant $C > 0$ so that the left side is less than or
equal to $C$ times the right side.  It is very easy to write down
explicit formulae for ``quasinorms'' which satisfy conditions (a) --
(d) and this weaker version of (e), and in fact this weaker version of
(e) is implied by the other conditions.  Also, any two quasinorms are
comparable, which is to say that each is bounded by a constant
multiple of the other.  Indeed, because of the homogeneity condition
(d), this statement can be reduced to one on a compact set not
containing the origin, where it follows from the continuity and
positivity of the quasinorms.

	Actual norms can be written down explicitly through simple but
carefully-chosen formulae, as in \cite{KR1}.  Another aspect of this
will be mentioned in a moment, but first let us define the distance
function associated to a norm or quasinorm.

	If $N$ is a norm or quasinorm on $H_n$, then we can define
an associated distance function $d_N(\cdot, \cdot)$ on $H_n$ by
\begin{equation}
	d_N\bigl((w,s), (z,t)\bigr) = N\bigl((w,s)^{-1} \circ (z,t)\bigr).
\end{equation}
By construction, this distance function is automatically invariant
under left translations on $H_n$, i.e.,
\begin{equation}
	d_N\bigl((y,u) \circ (w,s), (y,u) \circ (z,t)\bigr)
		= d_N\bigl((w,s), (z,t)\bigr)
\end{equation}
for all $(y,u), (w,s), (z,t) \in H_n$, simply because
\begin{equation}
	\bigl((y,u) \circ (w,s)\bigr)^{-1} \circ \bigl((y,u) \circ (z,t)\bigr)
		= (w,s)^{-1} \circ (z,t).
\end{equation}
We also have that $d(\cdot, \cdot)$ is nonnegative, equal to $0$ when the
two points in $H_n$ are the same, and is positive otherwise, because of
the corresponding properties of $N$.  Similarly,
\begin{equation}
	d_N\bigl((w,s), (z,t)\bigr)  = d_N\bigl((z,t), (w,s)\bigr),
\end{equation}
because of the symmetry property $N\bigl((w,s)^{-1}\bigr) = N(w,s)$ of $N$,
and 
\begin{equation}
	d_N\bigl(\delta_r(w,s), \delta_r(z,t)\bigr) 
		= r \, d_N\bigl((w,s), (z,t)\bigr)
\end{equation}
by the homogeneity property of $N$.

	If $N$ is a norm, then (\ref{triangle inequality for norm on
H_n}) implies that $d_N$ satisfies the usual triangle inequality for
metrics.  If $N$ is a quasinorm, then $d_N$ satisfies the weaker
version for quasimetrics, in which the right side is multiplied by a
fixed positive constant.  Just as different quasinorms on $H_n$ are
comparable, the corresponding distance functions are too, i.e.,
they are each bounded by a constant times the other.

	A basic and remarkable feature of the Heisenberg groups with
this geometry is that they are happy fractals.  In fact one can define
the distance between two points in terms of the infimum of the lengths
of certain paths between the two points, where the family of paths and
the notion of length enjoy left-invariance and homogeneity properties
which lead to the same kind of properties for the distance function as
above.  This kind of distance function can also be shown to be
compatible with the Euclidean topology on $H_n$.  These features imply
that this distance function is of the form $d_N$ for some $N$ as
above.  The triangle inequality for the distance function is a
consequence of its definition, and this leads to the triangle
inequality for the corresponding $N$.  A key subtlety in this approach
is that there is a sufficiently-ample supply of curves used in the
definition of the distance to connect arbitrary points in $H_n$,
because the curves are required to satisfy nontrivial conditions on
the directions of their tangent vectors.

	Let us return to the setting of an arbitrary norm $N$ on $H_n$.
The triangle inequality can be rewritten as
\begin{eqnarray}
	N(w,s) & \le & N(z,t) + d_N\bigl((w,s), (z,t)\bigr), 	\\
	N(z,t) & \le & N(w,s) + d_N\bigl((w,s), (z,t)\bigr)	\nonumber
\end{eqnarray}
for all $(w,s), (z,t) \in H_n$.  Thus
\begin{equation}
	|N(w,s) - N(z,t)| \le d_N\bigl((w,s), (z,t)\bigr)
\end{equation}
for all $(w,s), (z,t) \in H_n$.

	For $(w,s) \ne 0$, define $\phi(w,s)$ by
\begin{equation}
\label{def of phi(w,s)}
	\phi(w,s) = \delta_{N(w,s)^{-1}}(w,s).
\end{equation}
Thus
\begin{equation}
\label{N(phi(w,s)) = 1}
	N\bigl(\phi(w,s)\bigr) = 1
\end{equation}
by definition.

	If $(w,s)$, $(z,t)$ are both nonzero elements of $H_n$, then
\begin{eqnarray}
\lefteqn{d_N\bigl(\phi(w,s), \phi(z,t)\bigr)}	\\
	   & & \le d_N\bigl(\phi(w,s), \delta_{N(w,s)^{-1}}(z,t)\bigr)
		+ d_N\bigl(\delta_{N(w,s)^{-1}}(z,t), \phi(z,t)\bigr).
							\nonumber 
\end{eqnarray}
The first term on the right can be rewritten as
\begin{equation}
	d_N\bigl(\delta_{N(w,s)^{-1}}(w,s), \delta_{N(w,s)^{-1}}(z,t)\bigr)
		= N(w,s)^{-1} \, d_N\bigl((w,s), (z,t)\bigr),
\end{equation}
which is reasonable and nice for our purposes.  The second term on the
right can be rewritten as
\begin{equation}
	d_N\bigl(\delta_{N(z,t) \, N(w,s)^{-1}}(\phi(z,t)), \phi(z,t)\bigr).
\end{equation}
Let us think of this as being of the form
\begin{equation}
\label{d_N(delta_r(y,u), (y,u))}
	d_N\bigl(\delta_r(y,u), (y,u)\bigr),
\end{equation}
where $r$ is a positive real number and $(y,u) \in H_n$ satisfies
$N(y,u) = 1$.  Of course this expression is equal to $0$ when $r = 1$,
and one can be interested in getting a bound for it in terms of $r - 1$.

	Unfortunately one does not get a bound for
(\ref{d_N(delta_r(y,u), (y,u))}) like $O(|r-1|)$ in general, but more
like $O(\sqrt{|r-1|})$ for $r$ reasonably close to $1$.  The bottom
line is that the retraction $\phi$ onto the unit sphere for $N$ is not
Lipschitz, even in a small neighborhood of the sphere.  

	To look at it another way, although the dilation mapping
$\delta_r$ is Lipschitz with constant $r$ with respect to $d_N$ on
$H_n$, it does not have good Lipschitz properties as a function of
$r$, except on a small set.  This is in contrast to the case of
Euclidean geometry, where dilation by $r$ is uniformly Lipschitz as a
function of $r$ on bounded subsets.

	A closely related point is that while there are curves of
finite length joining the origin in $H_n$ to arbitrary elements of
$H_n$, the trajectories of the dilations do not have this property.

	Certainly one can expect that it is more difficult to have
Lipschitz retractions in the Heisenberg group than in Euclidean
spaces, and this indicates that this is so even for relatively simple
cases.

	Another basic mapping to consider is
\begin{equation}
	\psi(w,s) = \delta_{N(w,s)^{-2}}(w,s),
\end{equation}
which takes $H_n$ minus the origin to itself.  This mapping is a
reflection about the unit sphere for $N$, i.e., $\psi(w,s) = (w,s)$
when $N(w,s) = 1$, $N\bigl(\psi(w,s)\bigr) = N(w,s)^{-1}$, and
$\psi(\psi(w,s)) = (w,s)$.  Unlike the Euclidean case, there is once
again trouble with the Lipschitz condition even on a small
neighborhood of the unit sphere for $N$.

\section{Some happy fractals from Helsinki}
\label{Some happy fractals from Helsinki}
\setcounter{equation}{0}

	There are clearly numerous variations for the type of
construction about to be reviewed.  We shall focus on a simple family
with a lot of self-similarity.

	Let $N$ be an odd integer greater than or equal to $5$, and
let $\Sigma_0$ denote the boundary of the unit cube in ${\bf R}^3$.
Thus $\Sigma_0$ consists of $6$ two-dimensional squares, each with
sidelength $1$.

	In the first stage of the construction, we subdivide each of
these $6$ squares into $N^2$ squares with sidelength $1/N$.  For each
of the original $6$ squares, we make a modification with the square of
size $1/N$ in the middle.  The ``middle'' makes sense because $N$ is
odd.  Specifically, we remove the middle squares, and replace each one
with the union of the other $5$ squares in the boundary of the cube
with one face the middle square in question and which lies outside the
unit cube with which we started.  The surface that results from
$\Sigma_0$ by making these modifications is denoted $\Sigma_1$.

	This procedure can also be described as follows.  Let $R_0$
denote the unit cube, so that $\Sigma_0 = \partial R_0$.  Now define
$R_1$ to be the union of $R_0$ and the $6$ cubes with sidelength $1/N$
whose interiors are outside $R_0$ and which have a face which is
a middle square of a face of $R_0$.  The surface $\Sigma_1$ is the
boundary of $R_1$.

	Using the decomposition of the boundary described in the first
step, we can think of $\Sigma_1$ as the union of a bunch of
two-dimensional squares of sidelength $1/N$.  Namely, there are $6
\cdot (N^2 - 1) + 6 \cdot 5$ such squares.  For each of these squares,
we apply the same procedure as before.  That is, we divide each square
into $N^2$ squares of sidelength $1/N$ times the sidelength of the
squares that we have, so that the new squares have sidelength $1/N^2$
in this second step.  For each of the squares from the first step, we
make modifications only at the middle smaller squares just described,
one middle small square for each square from the second step.  Each of
these middle small squares is removed and replaced with the union of
$5$ squares of the same sidelength which are in the boundary of the
cube with interior outside $R_1$ and with one face being the small
middle square in question.  The result is a surface $\Sigma_2$
consisting of a bunch of squares of sidelength $1/N^2$.  The condition
$N \ge 5$ is helpful for keeping the modifications at different
places from bumping into each other or getting too close to doing that.

	One can also describe this in terms of adding a bunch of cubes
of sidelength $1/N^2$ to $R_1$, each with a face which is a middle
square of a square from the first step, to get a new region $R_2$.
The surface $\Sigma_2$ is the boundary of $R_2$.

	This process can be repeated indefinitely to get regions
$R_j$ and surfaces $\Sigma_j = \partial R_j$ for all nonnegative
integers $j$.  In the limit we can take $R$ to be the union of
the $R_j$'s, and $\Sigma$ to be the boundary of $R$, which is the
same as the Hausdorff limit of the $\Sigma_j$'s.

	Of course this procedure is completely analogous to ones in
the plane for producing snowflake curves.  However, one does not get
\emph{snowballs} in the technical sense introduced by Pekka Koskela,
because there are a lot of curves of finite length.  Indeed, whenever
a square is introduced in the construction, its four boundary segments
are kept intact for all future stages, and hence in the limit.  One
can verify that $\Sigma$ is a happy fractal.

\section{More on Lipschitz functions}
\label{More on Lipschitz functions}
\setcounter{equation}{0}

	Let $(M, d(x,y))$ be a metric space.  Suppose that $f(x)$ is a
real or complex-valued function on $M$, and that $L$ is a nonnegative
real number.  We say that $f$ is \emph{$L$-Lipschitz} if
\begin{equation}
\label{def of L-Lipschitz}
	|f(x) - f(y)| \le L \, d(x,y)
\end{equation}
for all $x, y \in M$.  Thus $f$ is Lipschitz if it is $L$-Lipschitz
for some $L$.  If $f$ is Lipschitz, then we define $\|f\|_{\Lip }$ to
be the supremum of
\begin{equation}
	\frac{|f(x) - f(y)|}{d(x,y)}
\end{equation}
over all $x, y \in M$, where this ratio is replaced with $0$ when $x =
y$.  In other words, $f$ is $\|f\|_{\Lip }$-Lipschitz when $f$ is
Lipschitz, and $\|f\|_{\Lip }$ is the smallest choice of $L$ for which
$f$ is $L$-Lipschitz.  Note that $\|\cdot \|_{\Lip }$ is a seminorm,
so that
\begin{equation}
	\|a \, f + b \, g \|_{\Lip } 
		\le |a| \, \|f\|_{\Lip } + |b| \, \|g\|_{\Lip }
\end{equation}
for all constants $a$, $b$ and Lipschitz functions $f$, $g$ on $M$.
Also, $\|f\|_{\Lip } = 0$ if and only if $f$ is a constant function on $M$.

	If $f$ and $g$ are real-valued $L$-Lipschitz functions on $M$,
then the maximum and minimum of $f$, $g$, which are denoted
$\max(f,g)$ and $\min(f,g)$, are $L$-Lipschitz functions too.  Let us
check this for $\max(f,g)$.  It is enough to show that
\begin{equation}
	\max(f,g)(x) - \max(f,g)(y) \le L \, d(x,y)
\end{equation}
for all $x, y \in M$, since one can interchange the roles of $x$ and
$y$ to get a corresponding lower bound for $\max(f,g)(x) -
\max(f,g)(y)$.  Assume, for the sake of definiteness, that
$\max(f,g)(x) = f(x)$.  Then we have
\begin{eqnarray}
	\ \max(f,g)(x) = f(x) & \le & f(y) + L \, d(x,y) 	\\
			& \le & \max(f,g)(y) + L \, d(x,y),	\nonumber
\end{eqnarray}
which is what we wanted.

	Here is a generalization of this fact.

\beginlemma
\label{lemma about sup of a family of L-Lipschitz functions}
Let $\{f_\sigma\}_{\sigma \in A}$ be a family of real-valued functions
on $M$ which are all $L$-Lipschitz for some $L \ge 0$.  Assume also
that there is point $p$ in $M$ such that the set of real numbers
$\{f_\sigma(p) : \sigma \in A\}$ is bounded from above.  Then the set
$\{f_\sigma(x) : \sigma \in A\}$ is bounded from above for every $x$
in $M$ (but not uniformly in $x$ in general), and $\sup \{f_\sigma(x)
: \sigma \in A\}$ is an $L$-Lipschitz function on $M$.
\end{lemma}

	Indeed, because $f_\sigma$ is $L$-Lipschitz for all $\sigma$
in $A$, we have that
\begin{equation}
	f_\sigma(x) \le f_\sigma(y) + L \, d(x,y)
\end{equation}
for all $x$, $y$ in $M$.  Applying this to $y = p$, we see that
$\{f_\sigma(x) : \sigma \in A\}$ is bounded from above for every $x$,
because of the analogous property for $p$.  If $F(x) = \sup
\{f_\sigma(x) : \sigma \in A\}$, then
\begin{equation}
	F(x) \le F(y) + L \, d(x,y)
\end{equation}
for all $x$, $y$ in $M$, so that $F$ is $L$-Lipschitz on $M$.  

	For the record, let us write down the analogous statement
for infima of $L$-Lipschitz functions.

\beginlemma
\label{lemma about inf of a family of L-Lipschitz functions}
Let $\{f_\sigma\}_{\sigma \in A}$ be a family of real-valued functions
on $M$ which are all $L$-Lipschitz for some $L \ge 0$.  Assume also
that there is point $q$ in $M$ such that the set of real numbers
$\{f_\sigma(q) : \sigma \in A\}$ is bounded from below.  Then the set
$\{f_\sigma(x) : \sigma \in A\}$ is bounded from below for every $x$
in $M$, and $\inf \{f_\sigma(x) : \sigma \in A\}$ is an $L$-Lipschitz
function on $M$.
\end{lemma}

	For any point $w$ in $M$, $d(x,w)$ defines a $1$-Lipschitz
function of $x$ on $M$.  This can be shown using the triangle inequality.
Suppose now that $f(x)$ is an $L$-Lipschitz function on $M$.  For
each $w \in M$, define $f_w(x) = f(w) + L \, d(x,w)$.  The fact
that $f$ is $L$-Lipschitz implies that
\begin{equation}
	f(x) \le f_w(x) \quad\hbox{for all } x, w \in M.
\end{equation}
Of course $f_x(x) = x$, and hence
\begin{equation}
	f(x) = \inf \{f_w(x) : w \in M\}.
\end{equation}
Each function $f_w(x)$ is $L$-Lipschitz in $x$, since $d(x,w)$ is
$1$-Lipschitz in $w$.  

	Similarly, we can set $\widetilde{f}_w(x) = f(x) - L \, d(x,w)$,
and then we have that
\begin{equation}
	f(x) = \sup \{\widetilde{f}_w(x) : w \in M \},
\end{equation}
and that $\widetilde{f}_w(x)$ is an $L$-Lipschitz function of $x$ for
every $w$.

	Here is a variant of these themes.  Let $E$ be a nonempty
subset of $M$, and suppose that $f$ is a real-valued function on $E$
which is $L$-Lipschitz, so that
\begin{equation}
	|f(x) - f(y)| \le L \, d(x,y)
\end{equation}
for all $x$, $y$ in $M$.  For each $w$ in $E$, set $f_w(x) = f(x) + L
\, d(x,w)$ and $\widetilde{f}_w(x) = f(x) - L \, d(x,w)$.  Consider
\begin{equation}
	F(x) = \inf \{f_w(x) : w \in E \}, \quad
		\widetilde{F}(x) = \sup \{\widetilde{f}_w(x) : w \in E \},
\end{equation}
for $x$ in $M$.  For the same reasons as before, $F(x) =
\widetilde{F}(x) = f(x)$ when $x$ lies in $E$.  Using Lemmas
\ref{lemma about sup of a family of L-Lipschitz functions} and
\ref{lemma about inf of a family of L-Lipschitz functions},
one can check that $F$ and $\widetilde{F}$ are $L$-Lipschitz
real-valued functions on all of $M$, i.e., they are extensions
of $f$ from $E$ to $M$ with the same Lipschitz constant $L$.

	If $H(x)$ is any other real-valued function on $M$ which
agrees with $f$ on $E$ and is $L$-Lipschitz, then
\begin{equation}
	\widetilde{f}_w(x) \le H(x) \le f_w(x)
\end{equation}
for all $w$ in $E$ and $x$ in $M$, and hence
\begin{equation}
	\widetilde{F}(x) \le H(x) \le F(x)
\end{equation}
for all $x$ in $M$.

\beginremark
\label{def of dist(x, S), and it is 1-Lipschitz}
{\rm 
If $S$ is any nonempty subset of $M$, define $\dist(x, S)$ for
$x$ in $M$ by
\begin{equation}
	\dist(x, S) = \inf_{y \in S} d(x,y).
\end{equation}
This function is always $1$-Lipschitz in $x$, as in Lemma \ref{lemma
about inf of a family of L-Lipschitz functions}.  }
\end{remark}

\section{Lipschitz functions of order $\alpha$}
\label{Lipschitz functions of order alpha}
\setcounter{equation}{0}

	Let $(M, d(x,y))$ be a metric space, and let $\alpha$ be a
positive real number.  A real or complex-valued function $f$ on $M$
is said to be \emph{Lipschitz of order $\alpha$} if there is nonnegative
real number $L$ such that
\begin{equation}
\label{def of Lip of order alpha}
	|f(x) - f(y)| \le L \, d(x,y)^\alpha
\end{equation}
for all $x, y \in M$.  This reduces to the Lipschitz condition
discussed in Section \ref{More on Lipschitz functions} when
$\alpha = 1$.  We shall sometimes write $\Lip \alpha$ for the
collection of Lipschitz functions of order $\alpha$, which might be
real or complex valued, depending on the context.  One also sometimes
refers to these functions as being ``H\"older continuous of order
$\alpha$''.

	If $f$ is Lipschitz of order $\alpha$, then we define
$\|f\|_{\Lip \alpha}$ to be the supremum of
\begin{equation}
	\frac{|f(x) - f(y)|}{d(x,y)^\alpha}
\end{equation}
over all $x, y \in M$, where this quantity is replaced with $0$ when
$x = y$.  In other words, $\|f\|_{\Lip \alpha}$ is the smallest choice
of $L$ so that (\ref{def of Lip of order alpha}) holds for all $x, y
\in M$.  This defines a seminorm on the space of Lipschitz functions
of order $\alpha$, as before, with $\|f\|_{\Lip \alpha} = 0$ if and
only if $f$ is constant.  Of course $\|f\|_{\Lip 1}$ is the same as
$\|f\|_{\Lip }$ from Section \ref{More on Lipschitz functions}.

	If $f$ and $g$ are real-valued functions on $M$ which are
Lipschitz of order $\alpha$ with constant $L$, then $\max(f,g)$ and
$\min(f,g)$ are also Lipschitz of order $\alpha$ with constant $L$.
This can be shown in the same manner as for $\alpha = 1$.  Similarly,
the analogues of Lemmas \ref{lemma about sup of a family of
L-Lipschitz functions} and \ref{lemma about inf of a family of
L-Lipschitz functions} for Lipschitz functions of order $\alpha$
hold for essentially the same reasons as before.

	However, if $\alpha > 1$, it may be that the only functions
that are Lipschitz of order $\alpha$ are the constant functions.  This
is the case when $M = {\bf R}^n$, for instance, equipped with the
standard Euclidean metric, because a function in $\Lip \alpha$ with
$\alpha > 1$ has first derivatives equal to $0$ everywhere.  Instead
of using derivatives, it is not hard to show that the function has to
be constant through more direct calculation too.  On any metric
space $M$, a function which is Lipschitz or order $\alpha$ with
$\alpha > 1$ is constant on every path of finite length.

	This problem does not occur when $\alpha < 1$.

\beginlemma
\label{if 0 < alpha le 1, then (a+b)^alpha le a^alpha + b^alpha}
If $0 < \alpha \le 1$ and $a$, $b$ are nonnegative real numbers, then
$(a+b)^\alpha \le a^\alpha + b^\alpha$.
\end{lemma}

	To see this, observe that 
\begin{equation}
	\max(a,b) \le (a^\alpha + b^\alpha)^{1/\alpha},
\end{equation}
and hence
\begin{eqnarray}
	a + b & \le & \max(a,b)^{1-\alpha} \, (a^\alpha + b^\alpha)	\\
		& \le & (a^\alpha + b^\alpha)^{1 + (1-\alpha)/\alpha}
		= (a^\alpha + b^\alpha)^{1/\alpha}.		\nonumber
\end{eqnarray}

\begincorollary
\label{d(x,y)^alpha a metric is d(x,y) is, 0 < alpha le 1}
If $(M, d(x,y))$ is a metric space and $\alpha$ is a real number
such that $0 < \alpha \le 1$, then $d(x,y)^\alpha$ also defines
a metric on $M$.
\end{corollary}

	This is easy to check.  The main point is that $d(x,y)^\alpha$
satisfies the triangle inequality, because of Lemma \ref{if 0 < alpha
le 1, then (a+b)^alpha le a^alpha + b^alpha} and the triangle
inequality for $d(x,y)$.

	A function $f$ on $M$ is Lipschitz of order $\alpha$ with
respect to the original metric $d(x,y)$ if and only if it is Lipschitz
of order $1$ with respect to $d(x,y)^\alpha$, and with the same norm.
In particular, for each $w$ in $M$, $d(x,w)^\alpha$ satisfies
(\ref{def of Lip of order alpha}) with $L = 1$ when $0 < \alpha \le
1$, because of the triangle inequality for $d(u,v)^\alpha$.

\section{Some functions on the real line}
\label{Some functions on the real line}
\setcounter{equation}{0}

	Fix $\alpha$, $0 < \alpha \le 1$.  For each nonnegative
integer $n$, consider the function
\begin{equation}
\label{2^{-n alpha} exp(2^n i x)}
	2^{-n \alpha} \exp(2^n \, i \, x)
\end{equation}
on the real line ${\bf R}$, where $\exp u$ denotes the usual exponential
$e^u$.  Let us estimate the $\Lip \alpha$ norm of this function.

	Recall that
\begin{equation}
\label{|exp(i u) - exp(i v)| le |u-v|}
	|\exp(i \, u) - \exp(i \, v)| \le |u-v|
\end{equation}
for all $u, v \in {\bf R}$.  Indeed, one can write $\exp(i \, u) -
\exp(i \, v)$ as the integral between $u$ and $v$ of the derivative of
$\exp(i \, t)$, and this derivative is $i \, \exp(i \, t)$, which has
modulus equal to $1$ at every point.

	Thus, for any $x, y \in {\bf R}$, we have that
\begin{equation}
  |2^{-n \alpha} \exp(2^n \, i \, x) - 2^{-n \alpha} \exp(2^n \, i \, y)|
	\le 2^{n (1-\alpha)} \, |x-y|.
\end{equation}
Of course
\begin{eqnarray}
\lefteqn{|2^{-n \alpha} \exp(2^n \, i \, x) 
   - 2^{-n \alpha} \exp(2^n \, i \, y)|} 	\\
	& & \le 
  2^{-n \alpha} |\exp(2^n \, i \, x)| + 2^{-n \alpha} |\exp(2^n \, i \, y)|
	= 2^{-n \alpha + 1}			\nonumber
\end{eqnarray}
as well.  As a result,
\begin{eqnarray}
\lefteqn{|2^{-n \alpha} \exp(2^n \, i \, x) 
   - 2^{-n \alpha} \exp(2^n \, i \, y)|} 	\\
	& & \le \Bigl(2^{n (1-\alpha)} \, |x-y| \Bigr)^{\alpha}
		\, \Bigl(2^{-n \alpha + 1} \Bigr)^{1-\alpha}
		= 2^{1-\alpha} \, |x-y|^\alpha.		\nonumber
\end{eqnarray}
This shows that the function (\ref{2^{-n alpha} exp(2^n i x)})
has $\Lip \alpha$ norm (with respect to the standard Euclidean
metric on ${\bf R}$) which is at most $2^{1-\alpha}$.  In the
opposite direction, if $2^n (x-y) = \pi$, then
\begin{eqnarray}
\lefteqn{|2^{-n \alpha} \exp(2^n \, i \, x) 
    - 2^{-n \alpha} \exp(2^n \, i \, y)|}		\\
	& & = 					
 2^{-n \alpha} |\exp(2^n \, i \, x)| + 2^{-n \alpha} |\exp(2^n \, i \, y)|
							\nonumber \\
	& & = 2^{-n \alpha + 1} = 2 \pi^{-\alpha} \, |x-y|^\alpha,
							\nonumber
\end{eqnarray}
so that the $\Lip \alpha$ norm is at least $2 \pi^{-\alpha}$.

	Now suppose that $f(x)$ is a complex-valued function on ${\bf
R}$ of the form
\begin{equation}
\label{f(x) = sum_{n=0}^infty a_n 2^{-n alpha} exp(2^n i x)}
	f(x) = \sum_{n=0}^\infty a_n \, 2^{-n \alpha} \exp(2^n \, i \, x),
\end{equation}
where the $a_n$'s are complex numbers.  We assume that the $a_n$'s are
bounded, which implies that the series defining $f(x)$ converges
absolutely for each $x$.  Set 
\begin{equation}
	A = \sup_{n \ge 0} |a_n|.
\end{equation}

	Let $m$ be a nonnegative integer.  For each $x$ in ${\bf R}$
we have that
\begin{equation}
\label{bound for sum from m to infty}
	\Bigl|\sum_{n=m}^\infty a_n \, 2^{-n \alpha} \exp(2^n \, i \, x) \Bigr|
	\le \sum_{n=m}^\infty A \, 2^{-n \alpha}
	= A \, (1 - 2^{-\alpha})^{-1} \, 2^{-m \alpha}.
\end{equation}
If $m \ge 1$ and $x, y \in {\bf R}$, then (\ref{|exp(i u) - exp(i v)|
le |u-v|}) yields
\begin{eqnarray}
\label{bound for difference of sums from 0 to m-1}
\lefteqn{\Bigl|\sum_{n=0}^{m-1} a_n \, 2^{-n \alpha} \exp(2^n \, i \, x) 
	- \sum_{n=0}^{m-1} a_n \, 2^{-n \alpha} \exp(2^n \, i \, y) \Bigr|} \\
   & & \le \sum_{n=0}^{m-1} A \, 2^{n (1-\alpha)} \, |x-y|
 \le A \, 2^{(m-1) (1-\alpha)} \Bigl(\sum_{j=0}^\infty 2^{-j (1-\alpha)}\Bigr)
		\, |x-y|				\nonumber \\	
   & & = A \, 2^{(m-1) (1-\alpha)} \, (1 - 2^{-(1-\alpha)})^{-1} \, |x-y|.
								\nonumber
\end{eqnarray}
Here we should assume that $\alpha < 1$, to get the convergence of
$\sum_{j=0}^\infty 2^{-j (1-\alpha)}$.  

	Fix $x, y \in {\bf R}$.  If $|x-y| > 1/2$, then we apply
(\ref{bound for sum from m to infty}) with $m = 0$ to both $x$ and $y$
to get that
\begin{eqnarray}
	|f(x) - f(y)| \le |f(x)| + |f(y)| 
		& \le & 2 \, A \, (1 - 2^{-\alpha})^{-1}		\\
	   & \le & 2^{1 + \alpha} \, A \, (1 - 2^{-\alpha})^{-1} |x-y|^\alpha.
								\nonumber
\end{eqnarray}
Assume now that $|x-y| \le 1/2$, and choose $m \in {\bf Z}_+$ so that
\begin{equation}
	2^{-m-1} < |x-y| \le 2^{-m}.
\end{equation}
Combining (\ref{bound for sum from m to infty}) and (\ref{bound for
difference of sums from 0 to m-1}), with (\ref{bound for sum from m to
infty}) applied to both $x$ and $y$, we obtain that
\begin{eqnarray}
\lefteqn{\quad |f(x) - f(y)|}						\\
  & & \le 2 \, A \, (1 - 2^{-\alpha})^{-1} \, 2^{-m \alpha}
	   + A \, 2^{(m-1) (1-\alpha)} \, (1 - 2^{-(1-\alpha)})^{-1} \, |x-y|
								\nonumber \\
  & & \le 2^{1+\alpha} \, A \, (1 - 2^{-\alpha})^{-1} \, |x-y|^\alpha
	+ A \, 2^{-(1-\alpha)} \, (1 - 2^{-(1-\alpha)})^{-1} \, |x-y|^\alpha.
								\nonumber
\end{eqnarray}
Therefore, for all $x, y \in {\bf R}$, we have that
\begin{eqnarray}
\lefteqn{|f(x) - f(y)|}						\\
  & & \le A (2^{1+\alpha} \, (1 - 2^{-\alpha})^{-1}
  + 2^{-(1-\alpha)} \, (1 - 2^{-(1-\alpha)})^{-1}) \, |x-y|^\alpha
								\nonumber
\end{eqnarray}
when $0 < \alpha < 1$.  In other words, $f$ is Lipschitz of order $\alpha$,
and
\begin{equation}
	\enspace \|f\|_{\Lip \alpha} 
		\le \Bigl(\sup_{n \ge 0} |a_n| \Bigr) \,
			(2^{1+\alpha} \, (1 - 2^{-\alpha})^{-1}
  			+ 2^{-(1-\alpha)} \, (1 - 2^{-(1-\alpha)})^{-1}).
\end{equation}

	To get an inequality going in the other direction we shall compute
as follows.  Let $\psi(x)$ be a function on ${\bf R}$ such that the
Fourier transform $\widehat{\psi}(\xi)$ of $\psi$,
\begin{equation}
   \widehat{\psi}(\xi) = \int_{\bf R} \exp (i \, \xi \, x ) \, \psi(x) \, dx
\end{equation}
is a smooth function which satisfies $\widehat{\psi}(1) = 1$ and
$\widehat{\psi}(\xi) = 0$ when $0 \le \xi \le 1/2$ and when $\xi \ge
2$.  One can do this with $\psi(x)$ in the Schwartz class of smooth
functions such that $\psi(x)$ and all of its derivatives are bounded
by constant multiples of $(1+|x|)^{-k}$ for every positive integer
$k$.

	For each nonnegative integer $j$, let us write $\psi_{2^j}(x)$
for the function $2^j \, \psi(2^j \, x)$.  Thus
\begin{equation}
	\widehat{\psi_{2^j}}(\xi) = \widehat{\psi}(2^{-j} \, \xi).
\end{equation}
In particular, $\widehat{\psi_{2^j}}(2^j) = 1$, and $\widehat{\psi_{2^j}}(2^l)
= 0$ when $l$ is a nonnegative integer different from $j$.  Hence
\begin{equation}
\label{int_{bf R} f(x) psi_{2^j}(x) dx = ....}
	\int_{\bf R} f(x) \, \psi_{2^j}(x) \, dx
   = \sum_{n=0}^\infty a_n \, 2^{-n \alpha} \, \widehat{\psi_{2^j}}(2^n)
	= a_j \, 2^{-j \alpha}.
\end{equation}
On the other hand,
\begin{equation}
	\int_{\bf R} \psi_{2^j}(x) \, dx = \widehat{\psi_{2^j}}(0)
		= \widehat{\psi}(0) = 0,
\end{equation}
so that
\begin{equation}
	\int_{\bf R} f(x) \, \psi_{2^j}(x) \, dx 
		= \int_{\bf R} (f(x) - f(0)) \, \psi_{2^j}(x) \, dx.
\end{equation}
Therefore
\begin{eqnarray}
	\Bigl|\int_{\bf R} f(x) \, \psi_{2^j}(x) \, dx \Bigr|
	& \le & \int_{\bf R} |f(x) - f(0)| \, |\psi_{2^j}(x)| \, dx	\\
  & \le & \|f\|_{\Lip \alpha} \int_{\bf R} |x|^\alpha \, |\psi_{2^j}(x)| \, dx
								\nonumber \\
  & = & \|f\|_{\Lip \alpha} \, 2^{-j \alpha} \, 
		\int_{\bf R} |x|^\alpha \, |\psi(x)| \, dx.	\nonumber
\end{eqnarray}
Combining this with (\ref{int_{bf R} f(x) psi_{2^j}(x) dx = ....}),
we obtain that
\begin{equation}
	|a_j| \le \|f\|_{\Lip \alpha} \, 
		\int_{\bf R} |x|^\alpha \, |\psi(x)| \, dx
\end{equation}
for all nonnegative integers $j$.  The integral on the right side
converges, because of the decay property of $\psi$.

	If $\alpha = 1$, then let us pass to the derivative and write
\begin{equation}
	f'(x) = \sum_{n=0}^\infty a_n \, i \, \exp{2^n \, i \, x}
\end{equation}
(where one should be careful about the meaning of $f'$ and of
this series).  This leads to
\begin{equation}
	\frac{1}{2\pi} \int_0^{2\pi} |f'(x)|^2 \, dx
		= \sum_{n=0}^\infty |a_n|^2.
\end{equation}
The main idea is that 
\begin{equation}
	\sum_{n=0}^\infty |a_n|^2 \le \|f\|_{\Lip 1}^2
\end{equation}
if $f$ is Lipschitz.  Conversely, if $\sum_{n=0}^\infty |a_n|^2 <
\infty$, then the derivative of $f$ exists in an $L^2$ sense, and
in fact one can show that $f'$ has ``vanishing mean oscillation''.

\section{Sums on general metric spaces}
\label{Sums on general metric spaces}
\setcounter{equation}{0}

	Let $(M, d(x,y))$ be a metric space.  For each integer $n$,
suppose that we have chosen a complex-valued Lipschitz function
$\beta_n(x)$ such that
\begin{equation}
	\sup_{x \in M} |\beta_n(x)| \le 1 \quad\hbox{and}\quad
		\|\beta\|_{\Lip } \le 2^n.
\end{equation}
Fix a real number $\alpha$, $0 < \alpha < 1$.

	Let $a_n$, $n \in {\bf Z}$ be a family (or doubly-infinite sequence)
of complex numbers which is bounded, and set 
\begin{equation}
	A = \sup_{n \in {\bf Z}} |a_n|.
\end{equation}
Consider 
\begin{equation}
	f(x) = \sum_{n \in {\bf Z}} a_n \, 2^{-n \alpha} \, \beta_n(x).
\end{equation}
The sum on the right side does not really converge in general, although
it would if we restricted ourselves to $n$ greater than any fixed number,
because of the bound on $\beta_n(x)$.  However, this sum does converge
``modulo constants'', in the sense that the sum in 
\begin{equation}
	f(x) - f(y) 
   = \sum_{n \in {\bf Z}} a_n \, 2^{-n \alpha} \, (\beta_n(x) - \beta_n(y)),
\end{equation}
converges absolutely for all $x$, $y$ in $M$.

	To see this, suppose that $k$ is any integer.  For $n \ge k$
we have that
\begin{equation}
	\sum_{n = k}^\infty |a_n| \, 2^{-n \alpha} \, |\beta_n(x)|
		\le A \, (1 - 2^{-\alpha})^{-1} \, 2^{-k \alpha},
\end{equation}
and similarly for $y$ instead of $x$.  For $n \le k-1$ we have that
\begin{eqnarray}
	\quad
   \sum_{n=-\infty}^{k-1} |a_n| \, 2^{-n \alpha} \, |\beta_n(x) - \beta_n(y)|
	& \le & A \sum_{n=-\infty}^{k-1} 2^{n (1-\alpha)} \, d(x,y)	\\
	& = & A \, 2^{(k-1)(1-\alpha)} \, (1 - 2^{-(1-\alpha)})^{-1} \, d(x,y).
								\nonumber
\end{eqnarray}
Thus
\begin{eqnarray}
\lefteqn{\sum_{n \in {\bf Z}} |a_n| \, 2^{-n \alpha} 
	   \, |\beta_n(x) - \beta_n(y)|}				\\
   & & \le A \, (1 - 2^{-\alpha})^{-1} \, 2^{-k \alpha}
     + A \, 2^{(k-1)(1-\alpha)} \, (1 - 2^{-(1-\alpha)})^{-1} \, d(x,y)
								\nonumber
\end{eqnarray}
for all $x, y \in M$ and $k \in {\bf Z}$.

\section{The Zygmund class on ${\bf R}$}
\label{The Zygmund class on R}
\setcounter{equation}{0}

	Let $f(x)$ be a real or complex-valued function on the real
line.  We say that $f$ lies in the \emph{Zygmund class} $Z$ if $f$ is
continuous and there is a nonnegative real number $L$ such that
\begin{equation}
\label{def of Zygmund condition}
	|f(x+h) + f(x-h) - 2 \, f(x)| \le L \, |h|
\end{equation}
for all $x, y \in {\bf R}$.  In this case, the seminorm $\|f\|_Z$ is
defined to be the supremum of
\begin{equation}
	\frac{|f(x+h) + f(x-h) - 2 \, f(x)|}{|h|}
\end{equation}
over all $x, h \in {\bf R}$ with $h \ne 0$.  This is the same as the
smallest $L$ so that (\ref{def of Zygmund condition}) holds.  Clearly
$f$ is in the Zygmund class when $f$ is Lipschitz (of order $1$), with
$\|f\|_Z \le 2 \, \|f\|_{\Lip }$.

	Suppose that $\{a_n\}_{n=0}^\infty$ is a bounded sequence of
complex numbers, and consider the function $f(x)$ on ${\bf R}$
defined by
\begin{equation}
	f(x) = \sum_{n=0}^\infty a_n \, 2^{-n} \, \exp(2^n \, i \, x).
\end{equation}
Let us check that $f$ lies in the Zygmund class, with $\|f\|_Z$ bounded
in terms of
\begin{equation}
	A = \sup_{n \ge 0} |a_n|.
\end{equation}
Note that $f$ is continuous.

	Observe that
\begin{eqnarray}
\lefteqn{|\exp(i (u+v)) + \exp(i (u-v)) - 2 \exp (i \, u)|}	\\
	& & = |\exp (i \, v) + \exp (-i \, v) - 2|		\nonumber
\end{eqnarray}
for all real numbers $u$, $v$, and that
\begin{equation}
	\exp (i \, v) + \exp (-i \, v) - 2 
		= \int_0^v i (\exp(i \, t) - \exp (-i \, t)) \, dt
\end{equation}
when $v \ge 0$.  Since $|\exp (i \, t) - \exp(-i \, t)| \le 2 \, t$ for
$t \ge 0$, we obtain that
\begin{equation}
	|\exp (i \, v) + \exp (-i \, v) - 2| \le \int_0^v 2 \, t \, dt = v^2.
\end{equation}
Hence
\begin{equation}
\label{|exp(i(u+v))+exp(i(u-v))-2 exp(i u)| le v^2}
	|\exp(i (u+v)) + \exp(i (u-v)) - 2 \exp (i \, u)|
		\le v^2,
\end{equation}
and this works for all real numbers $u$, $v$, since there is no real
difference between $v \ge 0$ and $v \le 0$.

	Let $x$ and $h$ be real numbers, and let $m$ be a nonnegative
integer.  From (\ref{|exp(i(u+v))+exp(i(u-v))-2 exp(i u)| le v^2})
we get that
\begin{eqnarray}
\lefteqn{\Bigl| \sum_{n=0}^m a_n \, 2^{-n} \, 
 (\exp(2^n \, i (x+h)) + \exp(2^n \, i (x-h)) - 2 \exp(2^n \, i \, x)) \Bigr|}
									\\
   & & \qquad\qquad\qquad\qquad
   \le A \sum_{n=0}^m 2^{-n} \, 2^{2n} \, |h|^2 \le A \, 2^{m+1} \, |h|^2.
								\nonumber
\end{eqnarray}
If $|h| \ge 1/2$, then 
\begin{eqnarray}
\lefteqn{|f(x+h) + f(x-h) - 2 \, f(x)|}			\\
  & &  \le |f(x+h)| + |f(x-h)| + 2 \, |f(x)| \le 4 \, A \le 8 \, A \, |h|.
							\nonumber
\end{eqnarray}
If $|h| \le 1/2$, then choose a positive integer $m$ such that
$2^{-m-1} \le |h| \le 2^{-m}$.  We can write $f(x+h) + f(x-h) - 2 \, f(x)$
as
\begin{eqnarray}
\lefteqn{\qquad \sum_{n=0}^m a_n \, 2^{-n} \, 
 (\exp(2^n \, i (x+h)) + \exp(2^n \, i (x-h)) - 2 \exp(2^n \, i \, x)) }
									\\
  & & + \sum_{n=m+1}^\infty a_n \, 2^{-n} \, 
 (\exp(2^n \, i (x+h)) + \exp(2^n \, i (x-h)) - 2 \exp(2^n \, i \, x)).
								\nonumber
\end{eqnarray}
This leads to
\begin{eqnarray}
\lefteqn{|f(x+h) + f(x-h) - 2 \, f(x)|}			\\
  & &  \le |f(x+h)| + |f(x-h)| + 2 \, |f(x)|		\nonumber \\
  & &  \le A \, 2^{m+1} \, |h|^2 + 4 \, A \, 2^{-m}	\nonumber \\
  & &  \le A \cdot 2 \cdot |h| + 4 \cdot A \cdot 2 \cdot |h| = 10 \, A \, |h|.
							\nonumber
\end{eqnarray}
This shows that $f$ lies in the Zygmund class, with constant less than
or equal to $10 \, A$.

\section{Approximation operators, 1}
\label{Approximation operators, 1}
\setcounter{equation}{0}

	Let $(M, d(x,y))$ be a metric space.  Fix a real number
$\alpha$, $0 < \alpha < 1$, and let $f$ be a real-valued function
on $M$ which is Lipschitz of order $\alpha$.  For each positive
real number $L$, define $A_L(f)$ by
\begin{equation}
\label{def of A_L(f)}
	A_L(f)(x) = \inf \{f(w) + L \, d(x,w) : w \in M \}
\end{equation}
for all $x$ in $M$.

	For arbitrary $x$, $w$ in $M$ we have that
\begin{equation}
\label{f(w) ge f(x) - ||f||_{Lip alpha} d(x,w)^alpha}
	f(w) \ge f(x) - \|f\|_{\Lip \alpha} \, d(x,w)^\alpha.
\end{equation}
As a result,
\begin{equation}
	f(w) + L \, d(x,w) \ge f(x) 
\end{equation}
when $L \, d(x,w)^{1-\alpha} \ge \|f\|_{\Lip \alpha}$.  Thus we
can rewrite (\ref{def of A_L(f)}) as
\begin{eqnarray}
\label{def of A_L(f), 2}
\lefteqn{A_L(f)(x) = }					\\
	& & \inf \{f(w) + L \, d(x,w) : 
		w \in M, \ L \, d(x,w)^{1-\alpha} \le \|f\|_{\Lip \alpha} \},
							\nonumber
\end{eqnarray}
i.e., one gets the same infimum over this smaller range of $w$'s.  In
particular, the set of numbers whose infimum is under consideration is
bounded from below, so that the infimum is finite.

	Because we can take $w = x$ in the infimum, we automatically
have that
\begin{equation}
\label{A_L(f)(x) le f(x)}
	A_L(f)(x) \le f(x)
\end{equation}
for all $x$ in $M$.  In the other direction, (\ref{f(w) ge f(x) -
||f||_{Lip alpha} d(x,w)^alpha}) and (\ref{def of A_L(f), 2})
lead to
\begin{eqnarray}
	A_L(f)(x) & \ge & f(x) - \|f\|_{\Lip \alpha} \, 
	   \biggl(\frac{\|f\|_{\Lip \alpha}}{L} \biggr)^{\alpha/(1-\alpha)} \\
    & = & f(x) - \|f\|_{\Lip \alpha}^{1/(1-\alpha)} \, L^{-\alpha/(1-\alpha)}.
								\nonumber
\end{eqnarray}
We also have that $A_L(f)$ is $L$-Lipschitz on $M$, as in Lemma
\ref{lemma about inf of a family of L-Lipschitz functions}.

	Suppose that $h(x)$ is a real-valued function on $M$ which
is $L$-Lipschitz and satisfies $h(x) \le f(x)$ for all $x$ in $M$.
Then
\begin{equation}
	h(x) \le h(w) + L \, d(x,w) \le f(w) + L \, d(x,w)
\end{equation}
for all $x$, $w$ in $M$.  Hence 
\begin{equation}
	h(x) \le A_L(f)(x)
\end{equation}
for all $x$ in $M$.

	Similarly, one can consider
\begin{equation}
\label{def of B_L(f)}
	B_L(f)(x) = \sup \{f(w) - L \, d(x,w) : w \in M \},
\end{equation}
and show that
\begin{eqnarray}
\label{def of B_L(f), 2}
\lefteqn{B_L(f)(x) = }					\\
	& & \sup \{f(w) - L \, d(x,w) : 
		w \in M, \ L \, d(x,w)^{1-\alpha} \le \|f\|_{\Lip \alpha} \}.
							\nonumber
\end{eqnarray}
This makes it clear that the supremum is finite.  As before,
\begin{equation}
	f(x) \le B_L(f)(x) 
   \le f(x) + \|f\|_{\Lip \alpha}^{1/(1-\alpha)} \, L^{-\alpha/(1-\alpha)},
\end{equation}
and $B_L(f)$ is $L$-Lipschitz.  If $h(x)$ is a real-valued function
on $M$ which is $L$-Lipschitz and satisfies $f(x) \le h(x)$ for all $x$
in $M$, then
\begin{equation}
	B_L(f)(x) \le h(x)
\end{equation}
for all $x$ in $M$.

\section{Approximation operators, 2}
\label{Approximation operators, 2}
\setcounter{equation}{0}

	Let $(M, d(x,y))$ be a metric space, and let $\mu$ be a
positive Borel measure on $M$.  We shall assume that $\mu$ is a 
\emph{doubling} measure, which means that there is a positive
real number $C$ such that
\begin{equation}
	\mu(B(x,2r)) \le C \, \mu(B(x,r))
\end{equation}
for all $x$ in $M$ and positive real numbers $r$, and that the
$\mu$-measure of any open ball is positive and finite.

	Let $t$ be a positive real number.  Define a function
$p_t(x,y)$ on $M \times M$ by
\begin{eqnarray}
	p_t(x,y) & = & 1 - t^{-1} d(x,y) \quad\hbox{when } d(x,y) \le t	\\
		 & = & 0 	  \qquad\qquad\qquad\enspace
				       \hbox{when } d(x,y) > t,		
								\nonumber
\end{eqnarray}
and put
\begin{equation}
	\rho_t(x) = \int_M p_t(x,y) \, d\mu(y).
\end{equation}
This is positive for every $x$ in $M$, because of the properties of $\mu$.
Also put
\begin{equation}
	\phi_t(x,y) = \rho_t(x)^{-1} \, p_t(x,y),
\end{equation}
so that 
\begin{equation}
\label{int_M phi_t(x,y) dmu(y) = 1}
	\int_M \phi_t(x,y) \, d\mu(y) = 1
\end{equation}
for all $x$ in $M$ by construction.

	Fix a real number $\alpha$, $0 < \alpha \le 1$, and let $f$
be a complex-valued function on $M$ which is Lipschitz of order $\alpha$.
Define $P_t(f)$ on $M$ by
\begin{equation}
	P_t(f)(x) = \int_M \phi_t(x,y) \, f(y) \, d\mu(y).
\end{equation}
Because of (\ref{int_M phi_t(x,y) dmu(y) = 1}), 
\begin{equation}
	P_t(f)(x) - f(x) = \int_M \phi_t(x,y) \, (f(y) - f(x)) \, d\mu(y),
\end{equation}
and hence
\begin{eqnarray}
  |P_t(f)(x) - f(x)| & \le & \int_M \phi_t(x,y) \, |f(y) - f(x)| \, d\mu(y)
									\\
      & \le & \int_M \phi_t(x,y) \, \|f\|_{\Lip \alpha} \, t^\alpha \, d\mu(y)
		= \|f\|_{\Lip \alpha} \, t^\alpha.		\nonumber
\end{eqnarray}
In the second step we employ the fact that $\phi_t(x,y) = 0$ when
$d(x,y) \ge t$.

	Suppose that $x$ and $z$ are elements of $M$, and consider
\begin{equation}
	|P_t(f)(x) - P_t(f)(z)|.
\end{equation}
If $d(x,z) \ge t$, then
\begin{eqnarray}
\lefteqn{|P_t(f)(x) - P_t(f)(z)|}				\\
   & & \le |P_t(f)(x) - f(x)| + |f(x) - f(z)| + |P_t(f)(z) - f(z)|	
								\nonumber \\
   & & \le  \|f\|_{\Lip \alpha} (2 \, t^\alpha + d(x,z)^\alpha)
		\le 3 \, t^{\alpha - 1} \, \|f\|_{\Lip \alpha} \, d(x,z).
								\nonumber
\end{eqnarray}
Assume instead that $d(x,z) \le t$.  In this case we write $P_t(f)(x)
- P_t(f)(z)$ as
\begin{eqnarray}
\lefteqn{\int_M (\phi_t(x,y) - \phi_t(z,y)) \, f(y) \, d\mu(y)}		\\
    & & = \int_M (\phi_t(x,y) - \phi_t(z,y)) \, (f(y) - f(x)) \, d\mu(y),
								\nonumber
\end{eqnarray}
using (\ref{int_M phi_t(x,y) dmu(y) = 1}).  This yields
\begin{eqnarray}
\lefteqn{|P_t(f)(x) - P_t(f)(z)|}					\\
    & & \le \int_M |\phi_t(x,y) - \phi_t(z,y)| \, |f(y) - f(x)| \, d\mu(y)
								\nonumber \\
    & & \le (2 t)^\alpha \, \|f\|_{\Lip \alpha} 
	\int_{\overline{B}(x, 2t)} |\phi_t(x,y) - \phi_t(z,y)| \, d\mu(y),
								\nonumber
\end{eqnarray}
where the second step relies on the observation that $\phi_t(x,y) -
\phi_t(z,y)$ is supported, as a function of $y$, in the set
\begin{equation}
	\overline{B}(x,t) \cup \overline{B}(z,t) \subseteq \overline{B}(x, 2t).
\end{equation}

	Of course
\begin{eqnarray}
\lefteqn{\quad \phi_t(x,y) - \phi_t(z,y)}				\\
    & & = (\rho_t(x)^{-1} - \rho_t(z)^{-1}) \, p_t(x,y)
		+ \rho_t(z)^{-1} \, (p_t(x,y) - p_t(z,y)).	\nonumber
\end{eqnarray}
Notice that
\begin{equation}
	|p_t(x,y) - p_t(z,y)| \le t^{-1} \, d(x,z)
\end{equation}
for all $y$ in $M$.  To see this, it is convenient to write $p_t(u,v)$
as $\lambda_t(d(u,v))$, where $\lambda_t(r)$ is defined for $r \ge 0$
by $\lambda_t(r) = 1 - t^{-1} \, r$ when $0 \le r \le t$, and
$\lambda_t(r) = 0$ when $r \ge t$.  It is easy to check that
$\lambda_t$ is $t^{-1}$-Lipschitz, and hence $\lambda_t(d(u,v))$ is
$t^{-1}$-Lipschitz on $M$ as a function of $u$ for each fixed $v$,
since $d(u,v)$ is $1$-Lipschitz as a function of $u$ for each fixed
$v$.  These computations and the doubling condition for $\mu$ permit
one to show that
\begin{equation}
	\int_{\overline{B}(x, 2t)} |\phi_t(x,y) - \phi_t(z,y)| \, d\mu(y)
		\le C_1 \, t^{-1} \, d(x,z)
\end{equation}
for some positive real number $C_1$ which does not depend on $x$, $z$,
or $t$.  (Exercise.)  Altogether, we obtain that
\begin{equation}
	\|P_t(f)\|_{\Lip 1} 
    \le \max(3, 2^\alpha \, C_1) \, t^{\alpha - 1} \, \|f\|_{\Lip \alpha}.
\end{equation}

\section{A kind of Calder\'on--Zygmund decomposition related to 
Lipschitz functions}
\label{A CZ decomposition related to Lipschitz functions}
\setcounter{equation}{0}

	Let $(M, d(x,y))$ be a metric space, and let $f$ be a
real-valued function on $M$.  Consider the associated maximal function
\begin{equation}
	N(f)(x) = \sup_{y \in M \atop y \ne x} \frac{|f(y) - f(x)|}{d(y,x)},
\end{equation}
where this supremum may be $+\infty$.

	Let $L$ be a positive real number, and put
\begin{equation}
	F_L = \{x \in M : N(f)(x) \le L \}.
\end{equation}
We shall assume for the rest of this section that
\begin{equation}
	F_L \ne \emptyset.
\end{equation}
As in Section \ref{Approximation operators, 1}, define $A_L(f)$
by
\begin{equation}
\label{def of A_L(f) repeated}
	A_L(f)(x) = \inf \{f(w) + L \, d(x,w) : w \in M \}.
\end{equation}
We shall address the finiteness of this infimum in a moment.
As before,
\begin{equation}
	A_L(f)(x) \le f(x)
\end{equation}
for all $x$ in $M$.

	If $u$ is any element of $F_L$, then
\begin{equation}
\label{|f(y) - f(u)| le L d(y,u), u in F_L}
	|f(y) - f(u)| \le L \, d(y,u)
\end{equation}
for all $y$ in $M$.  Let $x$ and $w$ be arbitrary points in $M$.
The preceding inequality implies that
\begin{equation}
	f(u) \le f(w) + L \, d(u,w),
\end{equation}
and hence
\begin{eqnarray}
	f(u) - L \, d(x,u) & \le & f(w) + L \, (d(u,w) - d(x,u))	\\
		& \le & f(w) + L \, d(x,w),			\nonumber
\end{eqnarray}
by the triangle inequality.  This yields
\begin{equation}
\label{f(u) - L d(x,u) le A_L(f)(x)}
	f(u) - L \, d(x,u) \le A_L(f)(x),
\end{equation}
which includes the finiteness of $A_L(f)(x)$.  If we take $x = u$,
then we get $f(u) \le A_L(f)(u)$, so that
\begin{equation}
	f(u) = A_L(f)(u) \quad\hbox{for all } u \in F_L.
\end{equation}
For $x \not\in F_L$, we obtain
\begin{equation}
	f(x) - 2 L \, d(x,u) \le A_L(f)(x)
\end{equation}
for all $u$ in $F_L$, by combining (\ref{f(u) - L d(x,u) le
A_L(f)(x)}) and (\ref{|f(y) - f(u)| le L d(y,u), u in F_L}) with $y =
x$.  In other words,
\begin{equation}
	f(x) - A_L(f)(x) \le 2 L \dist(x, F_L).
\end{equation}
Note that $A_L(f)$ is $L$-Lipschitz on $M$, by Lemma \ref{lemma about
inf of a family of L-Lipschitz functions}.

	In the same way, if
\begin{equation}
\label{def of B_L(f) repeated}
	B_L(f)(x) = \sup \{f(w) - L \, d(x,w) : w \in M \},
\end{equation}
then
\begin{equation}
	f(x) \le B_L(f)(x) \le f(x) + 2 L \dist(x,F_L)
\end{equation}
for all $x$ in $M$, and $B_L(f)$ is $L$-Lipschitz.

\section{A brief overview of ``atoms''}
\label{a brief overview of ``atoms''}
\setcounter{equation}{0}

	Let $(M, d(x,y))$ be a metric space, and let $s$ be a positive
real number.  We say that $(M, d(x,y))$ is \emph{Ahlfors-regular of
dimension $s$} if $M$ is complete as a metric space, and if there is a
positive Borel measure $\mu$ on $M$ such that
\begin{equation}
	C_1^{-1} \, r^s \le \mu(\overline{B}(x,r)) \le C_1 \, r^s
\end{equation}
for some positive real number $C_1$, all $x$ in $M$, and all $r > 0$
such that $r \le \diam M$ if $M$ is bounded.  

	As a basic example, if $M$ is $n$-dimensional Euclidean space
${\bf R}^n$ with the standard metric, and if $\mu$ is Lebesgue
measure, then in fact $\mu(\overline{B}(x,r))$ is equal to a constant
times $r^n$, where the constant is simply the volume of the unit ball.
More exotically, one can consider simply-connected nonabelian
nilpotent Lie groups, such as the Heisenberg groups.  For these spaces
one still has natural dilations as on Euclidean spaces, and Lebesgue
measure is compatible with both the group structure and the dilations,
in such a way that the measure of a ball of radius $r$ is equal to a
constant times $r^s$, where $s$ is now a geometric dimension that is
larger than the topological dimension.  Other examples include
fractals such as the Sierpinski gasket and carpet.

	Fix a metric space $(M, d(x,y))$ and a measure $\mu$ on $M$
satisfying the conditions in the definition of Ahlfors-regularity,
with dimension $s$.  The following fact is sometimes useful: there is
a constant $k_1 \ge 1$ so that if $x$ is an element of $M$ and $r$,
$R$ are positive numbers, with $r \le R$, then the ball
$\overline{B}(x,R)$ can be covered by a collection of at most $k_1 \,
(R/r)^s$ closed balls of radius $r$.  If $M$ is bounded, then we may
as well assume that $r < \diam M$ here, because $M$ is automatically
contained in a single ball with radius $\diam M$.  We may also assume
that $R \le \diam M$, since we could simply replace $R$ with $\diam M$
if $R$ is initially chosen to be larger than that.

	To establish the assertion in the preceding paragraph, let us
begin with a preliminary observation.  Suppose that $A$ is a subset of
$\overline{B}(x,R)$ such that $d(x,y) > r$ for all $x$, $y$ in $A$.
Then the number of elements of $A$ is at most $k_1 \, (R/r)^s$, if we
choose $k_1$ large enough (independently of $x$, $R$, and $r$).  Indeed,
\begin{equation}
	\sum_{a \in A} \mu(\overline{B}(a, r/2))
		= \mu\biggl(\bigcup_{a \in A} \overline{B}(a, r/2) \biggr)
		\le \mu(\overline{B}(x, 3R/2)),
\end{equation}
where the first equality uses the disjointness of the balls
$\overline{B}(a, r/2)$, $a \in A$.  The Ahlfors-regularity property
then applies to give a bound on the number of elements of $A$ of the
form $k_1 \, (R/r)^s$.  Now that we have such a bound, suppose that $A$
is also chosen so that the number of its elements is maximal.  Then
\begin{equation}
	\overline{B}(x,R) \subseteq \bigcup_{a \in A} \overline{B}(a,r).
\end{equation}
In other words, if $z$ is an element of $\overline{B}(x,R)$, then
$d(z,a) \le r$ for some $a$ in $A$, because otherwise we could add $z$
to $A$ to get a set which satisfies the same separation condition as
$A$, but which has $1$ more element.  This yields the original
assertion.

	In particular, closed and bounded subsets of $M$ are compact.
This uses the characterization of compactness in terms of completeness
and total boundedness, where the latter holds for bounded subsets of
$M$ by the result just discussed.

	Let us look at some special families of functions on $M$,
called \emph{atoms}, as in \cite{CW2}.  For the sake of definiteness,
we make the convention that a ``ball'' in $M$ means a closed ball
(with some center and radius), if nothing else is specified.  Suppose
that $p$ is a real number and $r$ is an extended real number such that
\begin{equation}
\label{conditions on p, r}
	0 < p \le 1, \ 1 \le r \le \infty, \ p < r.
\end{equation}
An integrable complex-valued function $a(x)$ on $M$ will be called a
\emph{$(p,r)$-atom} if it satisfies the following three conditions:
first, there is a ball $B$ in $M$ such that the support of $a$ is
contained in $B$, i.e., $a(x) = 0$ when $x \in M \backslash B$;
second,
\begin{equation}
\label{int_M a(x) d mu(x) = 0}
	\int_M a(x) \, d \mu(x) = 0;
\end{equation}
and third, 
\begin{equation}
\label{L^r average of a le mu(B)^{-1/p}}
	\biggl(\frac{1}{\mu(B)} \int_M |a(x)|^r \, d \mu(x) \biggr)^{1/r}
		\le \mu(B)^{-1/p}.
\end{equation}
If $r = \infty$, then (\ref{L^r average of a le mu(B)^{-1/p}})
is interpreted as meaning that the supremum (or essential supremum,
if one prefers) of $a$ is bounded by $\mu(B)^{-1/p}$.

	The size condition (\ref{L^r average of a le mu(B)^{-1/p}})
may seem a bit odd at first.  A basic point is that it implies
\begin{equation}
	\int_M |a(x)|^p \, d \mu(x) \le 1,
\end{equation}
by Jensen's inequality.  The index $r$ reflects a kind of regularity
of the atom, and notice that a $(p, r_1)$-atom is automatically a $(p,
r_2)$-atom when $r_1 \ge r_2$.  There are versions of this going in
the other direction, from $r_2$ to $r_1$, and we shall say more about
this soon.

	Suppose that $a(x)$ is a $(p,r)$-atom on $M$ and that
$\phi(x)$ lies in $\Lip \alpha$ on $M$ for some $\alpha$.  Consider
the integral
\begin{equation}
	\int_M a(x) \, \phi(x) \, d\mu(x).
\end{equation}
Let $B = \overline{B}(z,t)$ be the ball associated to $a(x)$ as
in the definition of an atom.  The preceding integral can be
written as
\begin{equation}
	\int_{\overline{B}(z,t)} a(x) \, (\phi(x) - \phi(z)) \, d\mu(x),
\end{equation}
using also (\ref{int_M a(x) d mu(x) = 0}).  Thus
\begin{eqnarray}
	\Bigl| \int_M a(x) \, \phi(x) \, d\mu(x) \Bigr|
   & \le & \int_{\overline{B}(z,t)} |a(x)| \, |\phi(x) - \phi(z)| \, d\mu(x)
									\\
 & \le & 
   \mu(\overline{B}(z,t))^{1 - (1/p)} \, t^\alpha \, \|\phi\|_{\Lip \alpha}.
								\nonumber
\end{eqnarray}
Ahlfors-regularity implies that
\begin{equation}
	\Bigl| \int_M a(x) \, \phi(x) \, d\mu(x) \Bigr|
  \le C_1^{1 - (1/p)} \, t^{(1 - (1/p)) s + \alpha} \, \|\phi\|_{\Lip \alpha}.
\end{equation}
In particular, 
\begin{equation}
\label{|int_M a(x) phi(x) dmu(x)| le C_1^{1 - (1/p)} ||phi||_{Lip alpha}}
	\Bigl| \int_M a(x) \, \phi(x) \, d\mu(x) \Bigr|
		  \le C_1^{1 - (1/p)} \, \|\phi\|_{\Lip \alpha}
\end{equation}
when $\alpha = ((1/p) - 1) \, s$.

	If we want to be able to choose $\alpha = ((1/p) - 1) \, s$
and have $\alpha \le 1$, then we are lead to the restriction
\begin{equation}
\label{p ge frac{s}{s+1}}
	p \ge \frac{s}{s+1}.
\end{equation}
Indeed, this condition does come up for some results, even if much of
the theory works without it.  There can also be some funny business at
the endpoint, so that one might wish to assume a strict inequality in
(\ref{p ge frac{s}{s+1}}), or some statements would have to be
modified when equality holds.

	In some situations this type of restriction is not really
necessary, perhaps with some adjustments.  Let us mention two basic
scenarios.  First, suppose that our metric space $M$ is something like
a self-similar Cantor set, such as the classical ``middle-thirds''
Cantor set.  In this case there are a lot of $\Lip \alpha$ functions
for all $\alpha > 0$, and, for that matter, there are a lot of
functions which are locally constant.  The computation giving
(\ref{|int_M a(x) phi(x) dmu(x)| le C_1^{1 - (1/p)} ||phi||_{Lip
alpha}}) still works when $\alpha > 1$, and this is true in general.

	On the other hand, if $M = {\bf R}^n$ with the standard
Euclidean metric, then there other ways to define classes of more
smooth functions, through conditions on higher derivatives.  In
connection with this, one can strengthen (\ref{int_M a(x) d mu(x) =
0}) by asking that the integral of an atom times a polynomial of
degree at most some number is equal to $0$.  If one does this, then
there are natural extensions of (\ref{|int_M a(x) phi(x) dmu(x)| le
C_1^{1 - (1/p)} ||phi||_{Lip alpha}}) for $\alpha > 1$, obtained by
subtracting a polynomial approximation to $\phi(x)$.  

	A basic manner in which atoms can be used is to test 
localization properties of linear operators.  Suppose that $T$
is a bounded linear operator on $L^2(M)$, and that $a$ is a
$(p,2)$-atom on $M$.  Consider
\begin{equation}
	T(a)
\end{equation}
(as well as $T^*(a)$, for that matter).  This is well-defined as an
element of $L^2(M)$, since $a$ lies in $L^2(M)$.  If $B =
\overline{B}(z,t)$ is the ball associated to $a$ in the definition of
an atom, then the estimate
\begin{eqnarray}
\lefteqn{\biggl(\frac{1}{\mu(B)} \int_M |T(a)(x)|^2 \, d \mu(x) \biggr)^{1/2}}
									\\
  & & \le  \|T\|_{2,2} \, 
	\biggl(\frac{1}{\mu(B)} \int_M |a(x)|^2 \, d \mu(x) \biggr)^{1/2}  
								\nonumber \\
	& & \le \|T\|_{2,2} \  \mu(B)^{-1/p}			\nonumber
\end{eqnarray}
provides about as much information about $T(a)$ around $B$, on $2 B =
\overline{B}(z,2t)$, say, as one might reasonably expect to have.
However, in many situations one can expect to have decay of $T(a)$
away from $B$, in such a way that
\begin{equation}
\label{||T(a)||_p le k}
	\|T(a)\|_p \le k 
\end{equation}
for some constant $k$ which does not depend on $a$.

	In this argument it is natural to take $r = 2$, but a basic
result in the theory is that one has some freedom to vary $r$.
Specifically, if $b$ is a $(p,r)$-atom on $M$, then it is possible
to write $b$ as
\begin{equation}
\label{b = sum_i beta_i b_i}
	b = \sum_i \beta_i \, b_i,
\end{equation}
where each $b_i$ is a $(p,\infty)$-atom, each $\beta_i$ is a complex
number, and $\sum_i |\beta_i|^p$ is bounded by a constant that does
not depend on $b$ (but which may depend on $p$ or $r$).  Let us give a
few hints about how one can approach this.  As an initial
approximation, one can try to write $b$ as
\begin{equation}
\label{b = beta' b' + sum_j gamma_j c_j}
	b = \beta' \, b' + \sum_j \gamma_j \, c_j,
\end{equation}
where $b'$ is a $(p,\infty)$-atom, $\beta'$ is a complex number such
that $|\beta'|$ is bounded by a constant that does not depend on $b$,
each $c_j$ is a $(p,r)$-atom, and $\sum_j |\gamma_j|^p \le 1/2$, say.
If one can do this, then one can repeat the process indefinitely to
get a decomposition as in (\ref{b = sum_i beta_i b_i}).  In order to
derive (\ref{b = beta' b' + sum_j gamma_j c_j}), the method of
Calder\'on--Zygmund decompositions can be employed.

	Recall that
\begin{equation}
	\Bigl(\sum_k \tau_k \Bigr)^p \le \sum_k \tau_k^p
\end{equation}
for nonnegative real numbers $\tau_k$ and $0 < p \le 1$.
As a consequence, if $\{f_k\}$ is a family of measurable
functions on $M$ such that 
\begin{equation}
	\int_M |f_k(x)|^p \, d\mu(x) \le 1 \qquad\hbox{for all } k,
\end{equation}
and if $\{\theta_k\}$ is a family of constants, then
\begin{equation}
	\int_M \Bigl|\sum_k \theta_k \, f_k(x) \Bigr|^p \, d\mu(x)
		\le \sum_k |\theta_k|^p.
\end{equation}
Because of this, bounds on $\sum_l |\alpha_l|^p$ are natural when
considering sums of the form $\sum_l \alpha_l \, a_l$, where the
$a_l$'s are $(p,r)$-atoms and the $\alpha_l$'s are constants.

	A fundamental theorem concerning atoms is the following.
Suppose that $T$ is a bounded linear operator on $L^2(M)$ again.  (One
could start as well with a bounded linear operator on some other $L^v$
space, with suitable adjustments.)  Suppose also that there is a
constant $k$ so that (\ref{||T(a)||_p le k}) holds for all
$(p,2)$-atoms, where $0 < p \le 1$, as before, or even simply for all
$(p,\infty)$-atoms.  Then $T$ determines a bounded linear operator
on $L^q$ for $1 < q < 2$.  This indicates how atoms are sufficiently
abundant to be useful.

	The proof of this theorem relies on an argument like the one
in Marcinkeiwicz interpolation.  In the traditional setting, one of
the main ingredients is to take a function $f$ in $L^q$ on $M$, and,
for a given positive real number $\lambda$, write it as $f_1 + f_2$,
where $f_1(x) = f(x)$ when $|f(x)| \le \lambda$, $f_1(x) = 0$ when
$|f(x)| > 0$, $f_2(x) = f(x)$ when $|f(x)| > \lambda$, and $f_2(x) =
0$ when $|f(x)| \le \lambda$.  Notice in particular that $f_1$ lies in
$L^w$ for all $w \ge q$, and that $f_2$ lies in $L^u$ for all $u \le
q$.  For the present purposes, the idea is to use decompositions which
are better behaved, with $f_2$ having a more precise form as a sum of
multiples of atoms.  The Calder\'on--Zygmund method is again
applicable, although it should be mentioned that one first works
with $(p,r)$-atoms with one choice of $r$, and then afterwards
makes a conversion to a larger $r$ using the results described before.

	In addition to considering the effect of $T$ on atoms,
one can consider the effect of $T^*$ on atoms, and this leads
to conclusions about $T$ on $L^q$ for $q > 2$, by duality.

\end{document}